\documentclass[10pt]{article}
\usepackage[latin1]{inputenc}
\usepackage{epsfig}
\usepackage{color}
\usepackage[british,english]{babel}
\usepackage{amsthm}
\usepackage{amsmath}
\usepackage{amsfonts}
\usepackage{amssymb}
\usepackage{graphicx}
\newtheorem{corollary}{Corollary}[section]

\newtheorem{lemma}[corollary]{Lemma}

\newtheorem{remark}[corollary]{Remark}
\newtheorem{theorem}[corollary]{Theorem}
\newfont{\sBlackboard}{msbm10 scaled 900}

\newcommand{\mylabel}[1]{\label{#1}
            \ifx\undefined\stillediting
            \else \fbox{$#1$}\fi }
\newcommand{\BE}{\begin{equation}}

\newcommand{\EEQ}{\end{equation}}
\newcommand{\rfb}[1]{\mbox{\rm
   (\ref{#1})}\ifx\undefined\stillediting\else:\fbox{$#1$}\fi}

\newfont{\Blackboard}{msbm10 scaled 1200}

\newfont{\roma}{cmr10 scaled 1200}

\def\CC{\rm \hbox{C\kern-.56em\raise.4ex
         \hbox{$\scriptscriptstyle |$}\kern+0.5 em }}

\def\n{|\kern -.05cm{|}\kern -.05cm{|}}

\def \noame{\noalign{\medskip}}
\newcommand{\mm}    {{\hbox{\hskip 0.5pt}}}

\newcommand{\bluff} {{\hbox{\raise 15pt \hbox{\mm}}}}

\newcommand{\ep}   {\varepsilon}
\usepackage{fancyhdr}

\makeatletter
\def\section{\@startsection {section}{1}{\z@}{-3.5ex plus -1ex minus
    -.2ex}{2.3ex plus .2ex}{\large\bf}}
\makeatother
\def\be{\begin{equation}}
\def\ee{\end{equation}}

\date{ }
\begin{document}
\thispagestyle{empty}
\title{\Large \bf Mathematical modeling of micropolar fluid flows\\ through a thin porous medium}\maketitle
\vspace{-2cm}
\begin{center}
Francisco Javier SU\'AREZ-GRAU\footnote{Departamento de Ecuaciones Diferenciales y An\'alisis Num\'erico. Facultad de Matem\'aticas. Universidad de Sevilla. 41012-Sevilla (Spain) fjsgrau@us.es}
 \end{center}


 \renewcommand{\abstractname} {\bf Abstract}
\begin{abstract} 
We study the flow of a micropolar fluid in a thin domain with microstructure, i.e. a thin domain with thickness $\varepsilon$ which is perforated by periodically distributed solid cylinders of size $a_\ep$. A main feature of this study is the dependence of the characteristic length of the micropolar fluid on the small parameters describing the geometry of the thin porous medium under consideration. Depending on the ratio of $a_\ep$ with respect to $\ep$, we derive three different generalized Darcy equations where the interaction between
the velocity and the microrotation fields is preserved.
\end{abstract}
\bigskip\noindent

\noindent {\small \bf AMS classification numbers:} 76A05, 76A20, 76M50, 76S05, 35B27, 35Q35.  \\
 

\section {Introduction}\label{S1}

Based on the micropolar fluid theory \cite{Eringen1, Eringen2}, which takes into account the effects of solid particles additive in a Newtonian fluid, we study flows of micropolar fluids in a thin domain which is perforated by periodically distributed
solid cylinders (microstructure) which is called thin porous medium ({\bf TPM}). This type of domains include two small parameters: one called $\ep$ is connected to the fluid film thickness and the other called $a_\ep$ to the microstructure representing the size of the cylinders and the interspatial distance between them.  The behavior of fluid  flows through {\bf TPM} has been studied extensively, mainly because of its importance in many industrial processes, see \cite{Fried, Jeon, Lund, Nord, Sinh, Tan}. However, the literature on non-Newtonian  micropolar fluid flows in this type of domains is far less complete, although these problems have now become of great practical relevance.  Therefore, the objetive of this paper is to derive generalized micropolar Darcy equations for the pressure depending on the magnitude of the parameters involving the {\bf TMP}.

For Newtonian fluids,  this problem has been addressed in  \cite{Fabricius} proving the existence of three types of {\bf TPM}:
\begin{itemize}
\item[-]  The proportionally thin porous medium ({\bf PTPM}), corresponding to the critical case when the cylinder height is proportional to the interspatial distance, with $\lambda$ the proportionality constant, that is $a_\ep\approx\ep$, with $a_\ep/\ep\to \lambda$, $0<\lambda<+\infty$. 
\item[-]  The homogeneously thin porous medium ({\bf HTPM}), corresponding  to the case when the cylinder height is much larger than interspatial distance, i.e. $a_\ep\ll \ep$ which is equivalent to $\lambda=0$.
\item[-]  The very thin porous medium ({\bf VTPM}), corresponding  to the case when the cylinder height is much smaller than the interspatial distance, i.e. $a_\ep\gg \ep$ which is equivalent to $\lambda=+\infty$.
\end{itemize}
In particular, starting from the Stokes system with body forces $f$, it can be  deduced that the flow is governed by a  2D Darcy equation for the pressure $p$ of the form
\begin{equation}\label{intro_Darcy}{\rm div}\left(K_\lambda(f-\nabla p)\right)=0\quad (0\leq \lambda\leq +\infty),\end{equation}
where $K_\lambda\in\mathbb{R}^{2\times 2}$ is a macroscopic quantity known as flow factor which takes into account the microstructure of the {\bf TPM}. Moreover, it holds that  the flow factor in the {\bf PTPM} is calculated by solving 3D Stokes local problems depending on the parameter $\lambda$, while  it is obtained in {\bf HTPM}  by solving 2D Stokes local problems and  by solving 2D local Hele-Shaw problems in the {\bf VTPM}, which it represents a considerable simplification.

These results were proved in  \cite{Fabricius} by using  the multiscale expansion method, which is a formal but powerful tool to analyze homogenization problems,  and later rigorously developed in \cite{AnguianoGrau2} by using an adaptation  of the periodic
unfolding method \cite{arbogast, Ciora, Ciora2}. This adaptation consists of a combination of the unfolding method in the horizontal variables  with a rescaling
in the height variable, in order to work with a domain of fixed height, and then to use suitable compactness results to pass
to the limit when the geometrical parameters $\ep$ and $a_\ep$ tend to zero. We remark that this adaptation was developed in \cite{AnguianoGrau}  to study the case of non-Newtonian power law fluids in the {\bf TPM} and  it was recently applied to the case of non-Newtonian Bingham fluids in  \cite{AnguianoBunoiu, AnguianoBunoiu2}.

Due to its importance in industrial and engineering applications, we consider a non-Newtonian micropolar  fluid  flow in {\bf TPM}  governed by  the linearized micropolar equations with body forces $f$ and body torque $g$.  By using the homogenization techniques developed in \cite{AnguianoGrau}, we derive that the flow is governed by a generalized 2D Darcy equation for the pressure $p$ of the form
\begin{equation}\label{intro_Darcy2}{\rm div}\left(K^{(1)}_\lambda(f-\nabla p)+ K^{(2)}_\lambda g\right)=0\quad (0\leq \lambda\leq +\infty),\end{equation}
where the flow factors $K^{(k)}\in\mathbb{R}^{2\times 2}$, $k=1,2$, are calculated depending on the cases as follows:  by solving 3D micropolar local problems depending on the parameter $\lambda$ in the {\bf PTPM}, by solving 2D micropolar local problems in the  {\bf HTPM} and by solving 2D local micropolar Reynolds problems in the case {\bf VTPM}.

As far as the author knows, this is the first attempt to carry out such an analysis for micropolar fluids, which could be instrumental for understanding the effects on the flows of  micropolar fluids and the microstructure of the domain. In view of that, more efficient numerical algorithms could be developed improving, hopefully, the known engineering practice.

The paper is organized as follows. In Section \ref{sec:notation}, we introduce some useful notation and in Section \ref{sec:estimates}, we give some a priori estimates for the dilated velocity, the microrotation and  the pressure, we introduce the extension of the unknowns to the whole domain $\Omega$, and finally we recall the version of the unfolding method necessary to pass to the limit in the next sections.  Namely, we analyze the case {\bf PTPM} in Section \ref{PTPM}, the  case {\bf HTPM} in Section \ref{HTPM} and the case {\bf VTPM} in Section \ref{VTPM}.

\section{Statement of the problem}\label{sec:notation}
In this section, we  first define the {\bf TPM} and some sets necessary to study the asymptotic behavior of the
solutions. Next, we introduce the problem considered and also the rescaled problem posed in a domain of  fixed height. We  finish this section giving the equivalent weak variational formulations.

\paragraph{ Definition of the TPM. }  A periodic porous medium is defined by a domain $\omega$ and an associated microstructure, or periodic cell  $Y'=(-1/2,1/2)^2$ which is made of two complementary parts: the fluid part $Y'_f$, and the solid part $Y'_s$ ($Y'_f\cup Y'_s=Y'$ and $Y'_f\cap Y'_s=\emptyset$). More precisely, we assume that $\omega$ is a smooth, bounded, connected set in $\mathbb{R}^2$ and that $Y'_s$ is an open connected subset of $Y'$ with a smooth boundary $\partial Y_s'$, such that $\overline {Y'_s}$ is strictly included in $Y'$. 

The microscale of a porous medium is given by a small positive number $a_\ep$. The domain $\omega$ is covered by a regular mesh of size $a_\ep$: for $k'\in \mathbb{Z}^2$, each cell $Y'_{k',a_\ep}=a_\ep k'+a_\ep Y'$ is divided in a fluid part $Y'_{f_{k'},a_\ep}$ and a solid part $Y'_{s_{k'},a_\ep}$, i.e. is similar to the unit cell $Y'$ rescaled to size $a_\ep$. We define $Y=Y'\times (0,1)\subset\mathbb{R}^3$, which is divided in a fluid part $Y_f$ and a solid part $Y_s$, and consequently $Y_{k',a_\ep}=Y'_{k',a_\ep}\times (0,1)\subset\mathbb{R}^3$, which is also divided in a fluid part $Y_{f_{k'},a_\ep}$, and a solid part $Y_{s_{k'},a_\ep}$.

We denote by $\tau(\overline Y'_{s_{k'},a_\ep})$ the set of all translated images of $\overline Y'_{s_{k'},a_\ep}$. The set $\tau(\overline Y'_{s_{k'},a_\ep})$ represents the solids in $\mathbb{R}^2$. The fluid part of the bottom $\omega_\ep\subset \mathbb{R}^2$ of the porous medium is defined by $\omega_\ep=\omega\setminus\bigcup_{k'\in\mathcal{K}_\ep}\overline Y_{s_{k'},a_\ep}$, where $\mathcal{K}_\ep=\{k'\in\mathbb{Z}^2\,:\, Y'_{k',a_\ep}\cap \omega\neq \emptyset\}$. The whole fluid part $\Omega_\ep\subset\mathbb{R}^3$  is defined by 
\begin{equation}\label{Omep}
\Omega_\ep=\{(x_1,x_2,x_3)\in\omega_\ep\times\mathbb{R}\,:\, 0<x_3<\ep\}.
\end{equation}
We make the assumption that the solids $\tau(\overline Y'_{s_{k'},a_\ep})$ do not intersect the boundary $\partial\omega$. We define $Y^\ep_{s_{k'},a_\ep}= Y'_{s_{k'},a_\ep}\times (0,\ep)$. Denote by $S_\ep$ the set of the solids contained in $\Omega_\ep$. Then, $S_\ep$ is a finite union of solids, i.e. $S_\ep=\bigcup_{k'\in\mathcal{K}_\ep}\overline Y^\ep_{s_{k'},a_\ep}$.

We define $\widetilde \Omega_\ep=\omega_\ep\times (0,1)$, $\Omega=\omega\times (0,1)$, and $Q_\ep=\omega\times (0,\ep)$. We observe that $\widetilde \Omega_\ep=\Omega\setminus \bigcup_{k'\in\mathcal{K}_\ep}\overline Y_{s_{k'},a_\ep}$, and we define $T_\ep=\bigcup_{k'\in\mathcal{K}_\ep}\overline Y_{s_{k'},a_\ep}$ as the set of the solids contained in $\widetilde\Omega_\ep$.

We remark that along this paper, the points $x\in \mathbb{R}^3$ will be decomposed as $x=(x',x_3)$ with $x'\in\mathbb{R}^2$, $x_3\in\mathbb{R}$. We also use the notation $x'$ to denote a generic vector of $\mathbb{R}^2$.

In order to apply the unfolding method, we will need the following notation. For $k'\in \mathbb{Z}^2$, we define $\kappa: \mathbb{R}^2\to \mathbb{Z}^2$ by
\begin{equation}\label{kappa_fun}
\kappa(x')=k' \Longleftrightarrow x'\in Y'_{k',1}\,.
\end{equation}
Remark that $\kappa$ is well defined up to a set of zero measure in $\mathbb{R}^2$ (the set $\cup_{k'\in\mathbb{Z}^2}\partial Y'_{k',1}$). Moreover, for every $a_\ep>0$, we have
$$\kappa\left({x'\over a_\ep}\right)=k' \Longleftrightarrow x'\in Y'_{k',a_\ep}\,.$$

The symbol $\#$ means $Y'$-periodicity. Thus, for example, the spaces $L^2_\#(Y)$, $H^1_\#(Y)$, $H^1_{0,\#}(Y)$ are defined by
$$L^2_\#(Y)=\left\{\varphi\in L^2_{\rm loc}(\mathbb{R}^2\times (0,1))\,:\,\int_Y|\varphi|^2dy<+\infty,\quad \varphi(y'+k',y_3)=\varphi(y),\quad \forall\, k'\in\mathbb{R}^2,\ \hbox{a.e. }y\in(\mathbb{R}^2\times (0,1)\right\},$$
$$H^1_\#(Y)=\left\{\varphi\in H^1_{\rm loc}(\mathbb{R}^2\times (0,1))\,:\, \varphi\in L^2_\#(Y),\ D\varphi\in L^2_\#(Y)^3\right\},$$
$$H^1_{0,\#}(Y)=\left\{\varphi\in H^1_\#(Y)\,:\,\varphi=0\quad\hbox{on }y_3=0,1\right\}.$$
We denote by $:$ the full contraction of two matrices; for $A=(a_{ij})_{1\leq i,j\leq 3}$ and $B=(b_{ij})_{1\leq i,j\leq 3}$, we have $A:B=\sum_{i,j=1}^3a_{ij}b_{ij}$.

We denote by $O_\ep$ a generic real sequence which tends to zero with $\ep$ and can change from line to line. We denote by $C$ a generic constant which can change from line to line.

\paragraph{The problem.}In $\Omega_\ep$  we consider the stationary flow of an incompressible micropolar fluid which is governed  by the following linearized micropolar system formulated in a non-dimensional form (see \cite{Luka} for more details)
\begin{equation}\label{system_1}
\left\{\begin{array}{rl}
-{\rm div}(D u_\ep)+\nabla p_\ep=2N^2{\rm rot}\,w_\ep+ f_\ep&\quad\hbox{in}\quad\Omega_\ep,\\
\noame
{\rm div}\,u_\ep=0&\quad\hbox{in}\quad\Omega_\ep,\\
\noame
-R_M{\rm div}(D w_\ep)+4N^2w_\ep=2N^2{\rm rot}\,u_\ep+g_\ep&\quad\hbox{in}\quad\Omega_\ep,
\end{array}\right.
\end{equation}
with homogeneous boundary conditions (it does not alter the generality of the problem under consideration),
\begin{equation}\label{bc_system_1}
u_\ep=w_\ep=0\quad \hbox{on}\quad\partial Q_\ep \cup \partial S_\ep\,.
\end{equation}

In system (\ref{system_1}), the velocity $u_\ep$, the pressure $p_\ep$ and the microrotation $w_\ep$ (i.e. the angular velocity field of rotation
 of particles) are unknown. Observe that equation of the linear momentum (\ref{system_1})$_1$ has the familiar form of the Stokes equation but it is coupled with equation of the angular momentum (\ref{system_1})$_3$,  which esentially describes the motion of the particles inside the microvolume as they undergo microrotational effects represented by the microrotational  vector $w_\ep$. Dimensionless (non-Newtonian) parameter $N^2$ with  $0\leq N\leq 1$  is called coupling number and it characterizes the coupling of the  linear and angular momentum equations.  When $N$ is identically zero, the equations  are decoupled and  equation of the linear momentum reduces to the classical Stokes equations which represent Newtonian fluids. The parameter $R_M$  is called characteristic length and it characterizes the interaction between the micropolar fluid and the microgeometry of the domain. It is small and has to be related to the other small geometrical parameters depending on the type of {\bf TPM}.

Under assumptions that $f_\ep,g_\ep\in L^2(\Omega_\ep)^3$, it is well known that problem (\ref{system_1})-(\ref{bc_system_1}) has a unique weak solution $(u_\ep,w_\ep,p_\ep)\in H^1_0(\Omega_\ep)^3\times H^1_0(\Omega_\ep)^3\times L^2_0(\Omega_\ep)$ (see  \cite{Luka}), where the space $L^2_0(\Omega_\ep)$ is the space of functions of $L^2(\Omega_\ep)$ with null integral.

Our aim is to study the asymptotic behavior of $u_\ep$, $w_\ep$ and $p_\ep$ when $\ep$ and $a_\ep$ tend to zero and identify homogenized models coupling the effects of the thickness of the domain and its microgeometry. For this purpose, as usual when we deal with thin domains, we use the dilatation in the variable $x_3$ given by
\begin{equation}\label{dilatacion}
y_3={x_3\over \ep}\,,
\end{equation}
in order to have the functions defined in the open set with fixed height   $\widetilde\Omega_\ep$.

Namely, we define  $\tilde u_\ep, \tilde w_\ep\in H^1_0(\widetilde\Omega_\ep)^3$ and $\tilde p_\ep\in L^2_0(\widetilde\Omega_\ep)$ by 
\begin{equation}\label{unk_dilat}
\tilde u_\ep(x',y_3)=u_\ep(x', \ep y_3),\quad \tilde w_\ep(x',y_3)=w_\ep(x',\ep y_3),\quad \tilde p_\ep(x',y_3)=p_\ep(x', \ep y_3),\quad\hbox{a.e. }(x',y_3)\in \widetilde\Omega_\ep\,.
\end{equation}
Let us introduce some notation which will be useful in the following. For a vectorial function $v=(v',v_3)$ and a scalar function $w$, we introduce the operators $D_{\ep}$, $\nabla_{ \ep}$ and ${\rm rot}_{ \ep}$ by 
\begin{eqnarray}
&\displaystyle(D_{ \ep}v)_{ij}=\partial_{x_j}v_i\hbox{ for }i=1,2,3,\ j=1,2,\quad  (D_{ \ep}v)_{i,3}={1\over  \ep}\partial_{y_3}v_i\hbox{ for }i=1,2,3,\nonumber&\\
\noame
&\displaystyle\nabla_{ \ep}w= \left(\nabla_{x'}w,{1\over  \ep}\partial_{y_3}w\right)^t,\quad {\rm div}_{ \ep}v={\rm div}_{x'}v'+{1\over  \ep}\partial_{y_3}v_3,\quad \displaystyle{\rm rot}_{ \ep}v=\left({\rm rot}_{x'}v_3+{1\over  \ep}{\rm rot}_{y_3}v',{\rm Rot}_{x'}v'\right)^t,\nonumber&\nonumber 
\end{eqnarray}
where, denoting $(v')^\perp=(-v_2,v_1)^t$, we define
\begin{equation}\label{def_rot}
{\rm rot}_{x'}v_3=(\partial_{x_2}v_3,-\partial_{x_1}v_3)^t,\quad {\rm rot}_{y_3}v'=(\partial_{y_3} v')^\perp,\quad {\rm Rot}_{x'}v'=\partial_{x_1}v_2-\partial_{x_2}v_1.
\end{equation}
Using the transformation (\ref{dilatacion}), the rescaled system (\ref{system_1})-(\ref{bc_system_1}) can be rewritten as
\begin{equation}\label{system_2}
\left\{
\begin{array}{rl}
-{\rm div}_{\ep}(D_{\ep} \tilde u_\ep)+\nabla_{\ep} \tilde p_\ep=2N^2{\rm rot}_{\ep}\,\tilde w_\ep+ \tilde f_\ep&\quad\hbox{in}\quad\widetilde\Omega_\ep,\\
\noame
{\rm div}_{\ep}\tilde u_\ep=0&\quad\hbox{in}\quad\widetilde \Omega_\ep,\\
\noame
-R_M{\rm div}_{\ep}(D_{\ep} \tilde w_\ep)+4N^2\tilde w_\ep=2N^2{\rm rot}_{\ep}\tilde u_\ep+\tilde g_\ep&\quad\hbox{in}\quad\widetilde\Omega_\ep\,,
\end{array}\right.
\end{equation}
with homogeneous boundary conditions
\begin{equation}\label{bc_system_2}
\tilde u_\ep=\tilde w_\ep=0\quad \hbox{on}\quad\partial \Omega \cup \partial T_\ep\,,
\end{equation}
where $\tilde f_\ep$ and $\tilde g_\ep$ are defined similarly as in (\ref{unk_dilat}). 

Our goal then is to describe the asymptotic behavior of this new sequences $\tilde u_\ep$, $\tilde w_\ep$ and $\tilde p_\ep$ when $\ep$ and $a_\ep$ tend to zero.  For this, it will be useful  to use the the equivalent weak variational formulation of  system (\ref{system_1})-(\ref{bc_system_1}) and  the rescaled system (\ref{system_2})-(\ref{bc_system_2}).
\paragraph{Weak variational formulations.}For  problem (\ref{system_1})-(\ref{bc_system_1}), the weak variational formulation is to find $u_\ep, w_\ep\in H^1_0(\Omega_\ep)^3$ and $p_\ep\in L^2_0(\Omega_\ep)$ such that
\begin{equation}\label{form_var_1}
\left\{\begin{array}{l}
\displaystyle \int_{\Omega_\ep}D u_\ep:D\varphi\,dx-\int_{\Omega_\ep}p_\ep\,{\rm div}\,\varphi\,dx=2N^2\int_{\Omega_\ep}{\rm rot}\,w_\ep\cdot \varphi\,dx+\int_{\Omega_\ep}f_\ep\cdot\varphi\,dx,\\
\noame
\displaystyle  R_M \int_{\Omega_\ep}D w_\ep:D\psi\,dx+4N^2\int_{\Omega_\ep}w_\ep\cdot\psi\,dx=2N^2\int_{\Omega_\ep}{\rm rot}\,u_\ep\cdot \psi\,dx+\int_{\Omega_\ep}g_\ep\cdot\psi\,dx\,,
\end{array}\right.
\end{equation}
for every $\varphi,\psi \in H^1_0(\Omega_\ep)^3$, and the equivalent  weak variational formulation for the rescaled system (\ref{system_2})-(\ref{bc_system_2}) is to find $\tilde u_\ep, \tilde w_\ep\in H^1_0(\widetilde \Omega_\ep)^3$ and $\tilde p_\ep\in L^2_0(\widetilde \Omega_\ep)$ such that

\begin{equation}\label{form_var_2}
\left\{\begin{array}{l}
\displaystyle \int_{\widetilde \Omega_\ep}D_{\ep} \tilde u_\ep:D_{\ep}\varphi\,dx'dy_3-\int_{\widetilde \Omega_\ep}\tilde p_\ep\,{\rm div}_{\ep}\varphi\,dx'dy_3=2N^2\int_{\widetilde \Omega_\ep}{\rm rot}_{\ep}\tilde w_\ep\cdot \varphi\,dx'dy_3+\int_{\widetilde \Omega_\ep}\tilde f_\ep\cdot\varphi\,dx'dy_3\,,\\
\noame
\displaystyle  R_M \int_{\widetilde \Omega_\ep}D_{\ep} \tilde w_\ep:D_{\ep}\psi\,dx'dy_3+4N^2\int_{\widetilde \Omega_\ep}\tilde w_\ep\cdot\psi\,dx'dy_3=2N^2\int_{\widetilde \Omega_\ep}{\rm rot}_{\ep}\tilde u_\ep\cdot \psi\,dx'dy_3+\int_{\widetilde \Omega_\ep}\tilde g_\ep\cdot\psi\,dx'dy_3\,,
\end{array}\right.
\end{equation}
for every $\varphi,\psi \in H^1_0(\widetilde \Omega_\ep)^3$.

\section{{\it A priori} estimates}\label{sec:estimates}
In the sequel we make the following assumptions concerning $f_\ep$, $g_\ep$, $R_M$ and $N$:
\begin{itemize}
\item[i)] in the cases  {\bf PTPM} and {\bf HTPM}, we assume
\begin{equation}\label{estim_f_g_cases_crit_sub}
f_\ep(x)=(f'(x'),0),\quad g_\ep(x)=(a_\ep g'(x'),0),\quad\hbox{a.e. }x\in\Omega_\ep,\quad \hbox{ where } f',g'\in L^2(\omega)^2,
\end{equation}
\begin{equation}\label{R_M}
N^2=\mathcal{O}(1),\quad R_M=a_\ep^2R_c\quad\hbox{with }R_c=\mathcal{O}(1)\,,
\end{equation}
\item[ii)] in the case {\bf VTPM}, we assume
\begin{equation}\label{estim_f_g_cases_sup}
f_\ep(x)=(f'(x'),0),\quad g_\ep(x)=(\varepsilon g'(x'),0),\quad\hbox{a.e. }x\in\Omega_\ep, \quad \hbox{ where } f',g'\in L^2(\omega)^2,
\end{equation}
\begin{equation}\label{R_M_super}
N^2=\mathcal{O}(1),\quad R_M=\ep^2R_c\quad\hbox{with }R_c=\mathcal{O}(1)\,.
\end{equation}
\end{itemize}
\begin{remark} We point out that in {\bf PTPM} and {\bf HTPM} the parameter $R_M$ is compared with  the size of the obstacles while in the case {\bf VTPM} with the film thickness, which are the most challenging ones and they answer to the question addressed in the paper, all preserve in the limit a strong coupling between velocity and microrotation. This choice is justified by many studies,  for example in the selected applications chapter in \cite{Luka} (see also \cite{BayadaChamGam, BayadaLuc}).

We also observe that due to the thickness of the domain, it is usual to assume that the vertical components of $f$ and $g$ can be neglected and, moreover they can be considered independent of the vertical variable. The parameters for $g_\ep$ are chosen to obtain appropriate estimates in each case.
\end{remark}

First, we recall the Poincar\'e inequality in a thin porous medium domain $\Omega_\varepsilon$ (see \cite{AnguianoGrau}).

\begin{lemma}\label{Poincare}
There exists a constant $C$ independent of $\varepsilon$, such that, 
\begin{itemize}
\item[i)] in the cases {\bf PTPM} and {\bf HTPM}, then 
\begin{equation}\label{p}
\left\Vert   v\right\Vert_{L^2( {\Omega}_{\varepsilon})^3}\leq Ca_{\varepsilon}\left\Vert D  v\right\Vert_{L^2( {\Omega}_{\varepsilon})^{{3\times3}}}, \quad \forall v\in H_0^1( {\Omega}_{\varepsilon})^3,
\end{equation}
\item[ii)] in the case {\bf VTPM}, then 
\begin{equation}\label{p_case2}
\left\Vert   v\right\Vert_{L^2({\Omega}_{\varepsilon})^3}\leq C{\varepsilon}\left\Vert D v\right\Vert_{L^2({\Omega}_{\varepsilon})^{{3\times3}}}, \quad \forall v\in H_0^1({\Omega}_{\varepsilon})^3.
\end{equation}
\end{itemize}
\end{lemma}

Next, we give the following results relating the derivative and the rotational (see \cite{DuvautLions}).
\begin{lemma}
 The following inequality holds
\begin{equation}\label{Gaffney}
\|{\rm rot}\,v\|_{L^2(\Omega_\ep)^{3}}\leq \|Dv\|_{L^2(\Omega_\ep)^{3\times 3}},\quad \forall v\in H_0^1( {\Omega}_{\varepsilon})^3.
\end{equation}
Moreover, if ${\rm div}\,v=0$ in $\Omega_\ep$, then it holds
\begin{equation}\label{Gaffney_div0}
\|{\rm rot}\,v\|_{L^2(\Omega_\ep)^{3}}=\|Dv\|_{L^2(\Omega_\ep)^{3\times 3}}.
\end{equation}
\end{lemma}

We start by obtaining some {\it a priori} estimates for $\tilde u_\ep$ and $\tilde w_\ep$.

\begin{lemma}\label{lemma_estimates}
There exists a constant $C$ independent of $\ep$, such that  the  rescaled solution $(\tilde u_\ep,\tilde w_\ep)$ of the problem (\ref{system_2})-(\ref{bc_system_2}) satisfies
\begin{itemize}
\item[i)] in the cases {\bf PTPM} and {\bf HTPM}, 
\begin{eqnarray}
&\|\tilde u_\ep\|_{L^2(\widetilde \Omega_\ep)^3}\leq Ca_\ep^{2},&\quad \|D_{\ep}\tilde u_\ep\|_{L^2(\widetilde \Omega_\ep)^{3\times 3}}\leq Ca_\ep\,,\label{estim_tilde_u_ep}\\
\noame
&\|\tilde w_\ep\|_{L^2(\widetilde \Omega_\ep)^3}\leq Ca_\varepsilon ,&\quad \|D_{\ep}\tilde w_\ep\|_{L^2(\widetilde \Omega_\ep)^{3\times 3}}\leq C \,.\label{estim_tilde_w_ep}
\end{eqnarray}
\item[ii)] in the case {\bf VTPM},
%
\begin{eqnarray}
&\|\tilde u_\ep\|_{L^2(\widetilde \Omega_\ep)^3}\leq C\ep^{2},&\quad \|D_{\ep}\tilde u_\ep\|_{L^2(\widetilde \Omega_\ep)^{3\times 3}}\leq C\ep\,,\label{estim_tilde_u_ep_super}\\
\noame
&\|\tilde w_\ep\|_{L^2(\widetilde \Omega_\ep)^3}\leq C\ep,&\quad \|D_{\ep}\tilde w_\ep\|_{L^2(\widetilde \Omega_\ep)^{3\times 3}}\leq C\,.\label{estim_tilde_w_ep_super}
\end{eqnarray}
\end{itemize}
\end{lemma}
{\bf Proof. } We analyze the different cases.
\begin{itemize} 
\item[i)] Cases {\bf PTPM} and {\bf HTPM}. We first obtain the estimates for the velocity. Taking $\varphi=u_\ep$ as test function in the first equation of (\ref{form_var_1}), taking into account $\int_{\Omega_\ep}{\rm rot}\,w_\ep\cdot u_\ep\,dx=\int_{\Omega_\ep}{\rm rot}\,u_\ep\cdot w_\ep\,dx$, applying Cauchy-Schwarz's inequality and from Lemma \ref{Poincare} and (\ref{Gaffney_div0}), we have 
\begin{eqnarray}
&\|Du_\ep\|^2_{L^2(\Omega_\ep)^{3\times 3}}& \displaystyle =2N^2\int_{\Omega_\ep}{\rm rot}\,w_\ep\cdot u_\ep\,dx+\int_{\Omega_\ep}f_\ep\cdot u_\ep\,dx\nonumber\\
\noame
&& \displaystyle =2N^2\int_{\Omega_\ep}w_\ep\cdot {\rm rot}\,u_\ep\,dx+\int_{\Omega_\ep}f'(x')\cdot u_\ep'\,dx\label{lem_estim_1}\\
\noame
&& \displaystyle \leq 2N^2 \|w_\ep\|_{L^2(\Omega_\ep)^3}\|D u_\ep\|_{L^2(\Omega_\ep)^{3\times 3}}+\ep^{1\over 2}a_\varepsilon  C\|f'\|_{L^2(\omega)^2}\|D u_\ep\|_{L^2(\Omega_\ep)^{3\times 3}}\,,\nonumber
\end{eqnarray}
which implies 
\begin{equation}\label{lem_estim_2}
\ep^{-{1\over 2}}a_\varepsilon^{-1}\|D u_\ep\|_{L^2(\Omega_\ep)^{3\times 3}}\leq \ep^{-{1\over 2}}a_\varepsilon^{-1}2N^2\|w_\ep\|_{L^2(\Omega_\ep)^3}+C\|f'\|_{L^2(\omega)^2}\,.
\end{equation}
Taking now $\psi=w_\ep$ as test function in the second equation of (\ref{form_var_1}), applying Cauchy-Schwarz's inequality and taking into account (\ref{estim_f_g_cases_crit_sub}) and (\ref{R_M}), we have
\begin{equation}\label{lem_estim_3}
\begin{array}{l}
\displaystyle
a_\ep^2R_c\|D w_\ep\|_{L^2(\Omega_\ep)^{3\times 3}}^2+4N^2\|w_\ep\|_{L^2(\Omega_\ep)^3}^2\\
\noame
\displaystyle\quad\quad  =2N^2\int_{\Omega_\ep}{\rm rot}\,u_\ep\cdot w_\ep\,dx+a_\ep\int_{\Omega_\ep}g'(x')\cdot w_\ep'\,dx\\
\noame
\displaystyle\quad\quad 
\leq 2N^2\|w_\ep\|_{L^2(\Omega_\ep)^3}\|D u_\ep\|_{L^2(\Omega_\ep)^{3\times 3}}+\ep^{1\over 2}a_\ep \|g'\|_{L^2(\omega)^2}\|w_\ep\|_{L^2(\Omega_\ep)^3}\,,
\end{array}
\end{equation}
which implies
\begin{equation}\label{lem_estim_4}
\ep^{-{1\over 2}}a_\varepsilon^{-1}2N^2\|w_\ep\|_{L^2(\Omega_\ep)^3}\leq \ep^{-{1\over 2}}a_\varepsilon^{-1}N^2\|D u_\ep\|_{L^2(\Omega_\ep)^{3\times 3}}+{1\over 2}\|g'\|_{L^2(\omega)^2}\,.
\end{equation}

Then, from (\ref{lem_estim_2}) and (\ref{lem_estim_4}), we conclude  that
$$\ep^{-{1\over 2}}a_\varepsilon^{-1}\|D u_\ep\|_{L^2(\Omega_\ep)^{3\times 3}}\leq {C\over 1-N^2}\|f'\|_{L^2(\omega)^2}+{1\over 2(1-N^2)}\|g'\|_{L^2(\omega)^2}\,,$$
which gives 
$$\|Du_\ep\|_{L^2(\Omega_\ep)^{3\times 3}}\leq Ca_\ep\ep^{1\over 2}.$$
 This together with  Lemma \ref{Poincare} gives 
 \begin{equation}\label{estim_u_ep}\|u_\ep\|_{L^2(\Omega_\ep)^3}\leq Ca_\ep^2\ep^{1\over 2},
 \end{equation}
and by means of the dilatation (\ref{dilatacion}) we get (\ref{estim_tilde_u_ep}).\\

Finally, we obtain the estimates for the microrotation. We use $\int_{\Omega_\ep}{\rm rot}\,u_\ep\cdot w_\ep\,dx=\int_{\Omega_\ep}{\rm rot}\,w_\ep\cdot u_\ep\,dx$  in (\ref{lem_estim_3}), Lemma \ref{Poincare} and  (\ref{Gaffney}), and proceeding as above we obtain 
\begin{equation}\label{lem_estim_5}
\begin{array}{l}
\displaystyle
a_\ep^2R_c\|D w_\ep\|_{L^2(\Omega_\ep)^{3\times 3}}^2+4N^2\|w_\ep\|_{L^2(\Omega_\ep)^3}^2\\
\noame
\displaystyle\quad\quad 
\leq 2N^2\|u_\ep\|_{L^2(\Omega_\ep)^3}\|D w_\ep\|_{L^2(\Omega_\ep)^{3\times 3}}+\ep^{1\over 2}a_\varepsilon^2C\|g'\|_{L^2(\omega)^2}\|Dw_\ep\|_{L^2(\Omega_\ep)^{3\times 3}}\,,
\end{array}
\end{equation}
which, by using the estimate of $u_\ep$ given in (\ref{estim_u_ep}), provides
$$a_\ep^2R_c\|D w_\ep\|_{L^2(\Omega_\ep)^{3\times 3}}\leq C\left(2N^2 \ep^{1\over 2}a_\varepsilon^2  + \ep^{1\over 2}a_\varepsilon^2\|g'\|_{L^2(\omega)^2}\right).$$
This implies 
$$\|w_\ep\|_{L^2(\Omega_\ep)^3}\leq Ca_\varepsilon\ep^{{1\over 2}},\quad \|Dw_\ep\|_{L^2(\Omega_\ep)^{3\times 3}}\leq C \ep^{{1\over 2}},$$ and by means of the dilatation, we get (\ref{estim_tilde_w_ep}).

\item[ii)] Case {\bf VTPM}. For the velocity, proceeding as above by taking into account (\ref{estim_f_g_cases_sup}) and (\ref{R_M_super}), we conclude
$$\ep^{-{3\over 2}}\|D u_\ep\|_{L^2(\Omega_\ep)^{3\times 3}}\leq \ep^{-{3\over 2}}2N^2\|w_\ep\|_{L^2(\Omega_\ep)^3}+C\|f'\|_{L^2(\omega)^2}\,,$$
and
$$\ep^{-{3\over 2}}2N^2\|w_\ep\|_{L^2(\Omega_\ep)^3}\leq \ep^{-{3\over 2}}N^2\|D u_\ep\|_{L^2(\Omega_\ep)^{3\times 3}}+{1\over 2}\|g'\|_{L^2(\omega)^2}\,.$$
Then, we deduce estimate 
$$\|u_\ep\|_{L^2(\Omega_\ep)^3}\leq C\ep^{5\over 2},\quad \|Du_\ep\|_{L^2(\Omega_\ep)^{3\times 3}}\leq C\ep^{3\over 2},$$
 and by means of the dilatation, we get (\ref{estim_tilde_u_ep_super}).

For the microrotation, similarly as the previous case, we get that 
$$\ep^2R_c\|D w_\ep\|_{L^2(\Omega_\ep)^{3\times 3}}\leq C\left(2N^2 \ep^{5\over 2}  + \ep^{5\over 2}\|g'\|_{L^2(\omega)^2}\right),$$
which implies 
$$\|w_\ep\|_{L^2(\Omega_\ep)^3}\leq C\ep^{{3\over 2}},\quad \|Dw_\ep\|_{L^2(\Omega_\ep)^{3\times 3}}\leq C\varepsilon^{{1\over 2}},$$ and by means of the dilatation we get (\ref{estim_tilde_w_ep_super}).
\end{itemize}
\par\hfill$\square$

\subsection{The extension of $(\tilde u_\ep,\tilde w_\ep,\tilde p_\ep)$ to the whole domain $\Omega$}
We extend the velocity $\tilde u_\varepsilon$ and the microrotation $\tilde w_\ep$ by zero to the $\Omega\setminus\widetilde \Omega_\ep$, and denote the extension by the same symbol. Obviously, estimates (\ref{estim_tilde_u_ep})-(\ref{estim_tilde_w_ep_super}) remain valid and   the extension of $\tilde u_\varepsilon$ is divergence free too.
 
 In order to extend the pressure to the whole domain $\Omega$, we use the mapping $R^\varepsilon$, defined in  Lemma 4.5 in \cite{AnguianoGrau} as $R^\ep_2$, which allows us to extend the pressure $p_\ep$ from $\Omega_\ep$ to $Q_\ep$ by introducing $F_\ep$ in $H^{-1}(Q_\ep)^3$ in the following way  (brackets are for the duality products between $H^{-1}$ and $H^1_0$):
\begin{equation}\label{F}\langle F_\varepsilon, \varphi\rangle_{Q_\varepsilon}=\langle \nabla p_\varepsilon, R^\varepsilon \varphi\rangle_{\Omega_\varepsilon}\,,\quad \hbox{for any }\varphi\in H^{1}_0(Q_\varepsilon)^3\,.
\end{equation}
Using Lemma \ref{lemma_estimates} for fixed $\varepsilon$ and $a_\ep$, we see that it i a bounded functional on $H^1_0(Q_\varepsilon)$ (see the proof of Lemma \ref{lemma_est_P} below), and in fact $F_\varepsilon\in H^{-1}(Q_\varepsilon)^3$. Moreover, ${\rm div} \varphi=0$ implies $\left\langle F_{\varepsilon},\varphi\right\rangle_{Q_\varepsilon}=0\,,$ and the DeRham theorem gives the existence of $P_\varepsilon$ in $L^{2}_0(Q_\varepsilon)$ with $F_\varepsilon=\nabla P_\varepsilon$.

We calcule the right hand side of (\ref{F}) by using the first equation of (\ref{form_var_1}) and we have
\begin{equation}\label{equality_duality}
\begin{array}{rl}
\displaystyle
\left\langle F_{\varepsilon},\varphi\right\rangle_{Q_\varepsilon}=&\displaystyle
-\int_{\Omega_\varepsilon}D u_\ep : D R^{\varepsilon}\varphi\,dx + 2N^2 \int_{\Omega_\ep}{\rm rot}\,w_\ep\cdot R^\ep \varphi\,dx+\int_{\Omega_\varepsilon} f'(x')\cdot (R^{\varepsilon}\varphi)'\,dx\,.
\end{array}\end{equation}
We get for any $\tilde \varphi\in H^1_0(\Omega)^3$ where $\tilde \varphi(x', y_3)=\varphi(x',\ep y_3)$, using the change of variables (\ref{dilatacion}), that
$$\begin{array}{rl}\displaystyle\langle \nabla_{\varepsilon}\tilde P_\varepsilon, \tilde\varphi\rangle_{\Omega}&\displaystyle
=-\int_{\Omega}\tilde P_\varepsilon\,{\rm div}_{\varepsilon}\,\tilde\varphi\,dx'dy_3
=-\varepsilon^{-1}\int_{Q_\varepsilon}P_\varepsilon\,{\rm div}\,\varphi\,dx=\varepsilon^{-1}\langle \nabla P_\varepsilon, \varphi\rangle_{Q_\varepsilon}\,.
\end{array}$$

Then, using the identification (\ref{equality_duality}) of $F_\varepsilon$, we have
$$\begin{array}{l}\displaystyle\langle \nabla_{\varepsilon}\tilde P_\varepsilon, \tilde\varphi\rangle_{\Omega}\displaystyle
=\varepsilon^{-1}\Big(-\int_{\Omega_\varepsilon}D u_\ep : D R^{\varepsilon} \varphi\,dx
+2N^2\int_{\Omega_\ep}{\rm rot}\,w_\ep\cdot R^\ep \varphi \,dx+\int_{\Omega_\varepsilon} f'(x')\cdot (R^{\varepsilon} \varphi)'\,dx\Big)\,,
\end{array}$$
and applying the change of variables (\ref{dilatacion}), we obtain 
\begin{equation}\label{extension_1}
\begin{array}{l}\displaystyle\langle \nabla_{\varepsilon}\tilde P_\varepsilon, \tilde\varphi\rangle_{\Omega}
=\displaystyle- \int_{\widetilde \Omega_\varepsilon}D_{\ep} \tilde u_\ep : D_{\ep} \tilde R^{\varepsilon} \tilde \varphi\,dx'dy_3
\displaystyle +2N^2\int_{\widetilde \Omega_\ep}{\rm rot}_{\ep}\tilde w_\ep\cdot \tilde R^\ep \tilde \varphi\,dx'dy_3+\int_{\widetilde \Omega_\varepsilon} f(x')\cdot (\tilde R^{\varepsilon} \tilde \varphi)'\,dx'dy_3\,,
\end{array}
\end{equation}
where $\tilde R^\ep\tilde \varphi=R^\ep \varphi$ for any $\tilde \varphi \in H^1_0(\Omega)^3$.

Now, we estimate the right-hand side of (\ref{extension_1}).

\begin{lemma}\label{lemma_est_P}
There exists a constant $C>0$ independent of $\ep$, such that the extension $\tilde P_\ep\in L^2_0(\Omega)$ of the pressure $\tilde p_\ep$ satisfies 
\begin{equation}\label{esti_P}
\| \tilde{P}_{\varepsilon} \|_{L_0^{2}(\Omega)}\leq {C}.
\end{equation}
\end{lemma}
{\bf Proof. } From the proof of Lemma 4.6-(i) in \cite{AnguianoGrau}, we have that $\tilde R^\ep\tilde\varphi$ satisfies the following estimates
$$\left\{\begin{array}{l}\displaystyle\|\tilde R^\ep(\tilde\varphi)\|^2_{L^2(\widetilde\Omega_\ep)^3}\leq C\left({1\over a_\ep^2}\|\tilde\varphi\|^2_{L^2(\Omega)^3} 
+ \|D_{x'}\tilde\varphi\|^2_{L^2(\Omega)^{3\times 2}}+  {1\over a_\ep^2}\|\partial_{y_3}\tilde\varphi\|^2_{L^2(\Omega)^{3}}\right)\leq C\|\tilde\varphi\|_{H^1_0(\Omega)^3},\\
\noame
\displaystyle
\|D_{x'}\tilde R^\ep\tilde\varphi\|^2_{L^2(\widetilde\Omega_\ep)^{3\times 2}}\leq C\left({1\over a_\ep^2}\|\tilde\varphi\|^2_{L^2(\Omega)^3} 
+ \|D_{x'}\tilde\varphi\|^2_{L^2(\Omega)^{3\times 2}}+  {1\over a_\ep^2}\|\partial_{y_3}\tilde\varphi\|^2_{L^2(\Omega)^{3}}\right),\\
\noame
\|\partial_{y_3}\tilde R^\ep\tilde\varphi\|^2_{L^2(\widetilde\Omega_\ep)^{3}}\leq C\left(\|\tilde\varphi\|^2_{L^2(\Omega)^3} 
+ a_\ep^2\|D_{x'}\tilde\varphi\|^2_{L^2(\Omega)^{2\times 3}}+\|\partial_{y_3}\tilde\varphi\|^2_{L^2(\Omega)^{3}}\right).
\end{array}\right.
$$
This implies that 

$$\|\tilde R^\ep(\tilde\varphi)\|_{L^2(\widetilde\Omega_\ep)^3}\leq C\left(\|\tilde\varphi\|^2_{L^2(\Omega)^3} 
+ a_\ep\|D_{x'}\tilde\varphi\|_{L^2(\Omega)^{3\times 2}}+ \|\partial_{y_3}\tilde\varphi\|^2_{L^2(\Omega)^{3}}\right)\leq C\|\tilde\varphi\|_{H^1_0(\Omega)^3}.$$
 Moreover, in the case  {\bf PTPM} and {\bf HTPM},
$$
\begin{array}{l}\displaystyle
\|D_{\ep}\tilde R^\ep\tilde\varphi\|_{L^2(\widetilde\Omega_\ep)^{3\times 3}}\leq C\left({1\over a_\ep}\|\tilde\varphi\|_{L^2(\Omega)^3} + \|D_{x'}\tilde\varphi\|_{L^2(\Omega)^{3\times 2}}+{1\over a_\ep}\|\partial_{y_3}\tilde \varphi\|_{L^2(\Omega)^3}\right)\leq {C\over a_\ep}\|\tilde\varphi\|_{H^1_0(\Omega)^3},
\end{array}
$$
 and in the case {\bf VTPM},
$$\|D_{\ep}\tilde R^\ep \tilde\varphi\|_{L^2(\widetilde\Omega_\ep)^{3\times 3}}\leq C\left({1\over \ep}\|\tilde\varphi\|_{L^2(\Omega)^3} + {a_\ep\over \ep}\|D_{x'}\tilde\varphi\|_{L^2(\Omega)^{3\times 3}}+ {1\over \ep}\|\partial_{y_3}\tilde \varphi\|_{L^2(\Omega)^3}\right)\leq {C\over \ep}\|\tilde\varphi\|_{H^1_0(\Omega)^3} $$

Thus, in the cases {\bf PTMP} and {\bf HTPM}, by using estimates for $D_{\ep}\tilde u_\ep$ in (\ref{estim_tilde_u_ep}), for $D_{\ep}w_\ep$ in (\ref{estim_tilde_w_ep}) and  $f'\in L^2(\omega)^2$,  we respectively obtain
$$
\begin{array}{l}
\displaystyle
\left|\int_{\widetilde\Omega_\ep}D_{\ep}\tilde u_\ep:D_{\ep}\tilde R^\ep\tilde\varphi\,dx'dy_3\right|\leq Ca_\ep\|D_{\ep}\tilde R^\ep\tilde\varphi\|_{L^2(\widetilde \Omega_\ep)^{3\times 3}}\leq C\|\tilde\varphi\|_{H^1_0(\Omega)^3},\\
\noame
\displaystyle \left|\int_{\widetilde\Omega_\ep}{\rm rot}_{\ep}w_\ep\cdot \tilde R^\ep \tilde\varphi \,dx'dy_3\right|\leq 
\|D_{\ep}\tilde w_\ep\|_{L^2(\widetilde\Omega_\ep)^{3\times 3}}\|\tilde R^\ep \tilde\varphi \|_{L^2(\widetilde\Omega_\ep)^3}\leq Ca_\ep^2\ep^{-2}\|\tilde R^\ep \tilde\varphi \|_{L^2(\widetilde\Omega_\ep)^3}\leq C\|\tilde\varphi\|_{H^1_0(\Omega)^3},\\
\noame\displaystyle
 \left|\int_{\widetilde\Omega_\ep}f'\cdot \tilde R^\ep \tilde\varphi \,dx'dy_3\right|\leq C\|\tilde R^\ep \tilde\varphi \|_{L^2(\widetilde\Omega_\ep)^3}\leq C\|\tilde\varphi\|_{H^1_0(\Omega)^3}\,,
\end{array}
$$
which together with (\ref{extension_1}) gives $\|\nabla_{\ep}\tilde P_\ep\|_{L^2(\Omega)^3}\leq C$. By using the Ne${\breve{\rm c}}$as inequality there exists a representative $\tilde P_\ep\in L^2_0(\Omega)$ such that
$$\|\tilde P_\ep\|_{L^2(\Omega)}\leq C\|\nabla\tilde P_\ep\|_{H^{-1}(\Omega)^3}\leq C\|\nabla_{\ep}\tilde P_\ep\|_{H^{-1}(\Omega)^3},$$
which implies (\ref{esti_P}).

 Finally, proceeding similarly in the case {\bf VTPM}, by using now estimates for $D_{\ep}\tilde u_\ep$ in (\ref{estim_tilde_u_ep_super}) and for $D_{\ep}w_\ep$ in (\ref{estim_tilde_w_ep_super}),  we also obtain  (\ref{esti_P}).
\par\hfill$\square$

\subsection{Adaptation of the unfolding method}\label{sec:unfolding}
The change of variables (\ref{dilatacion}) does not provide the information we need about the behavior of $\tilde u_\ep$ and $\tilde w_\ep$ in the microstructure associated to $\widetilde\Omega_\ep$. To solve this difficulty, we use an adaptation of the unfolding method (see \cite{arbogast, Ciora, Ciora2} for more details) introduced in \cite{AnguianoGrau}.
\\

Let us recall that this adaptation of the unfolding method divides the domain $\widetilde\Omega_\ep$ in cubes of lateral length $a_\ep$ and vertical length $1$. Thus, given $(\tilde{u}_{\varepsilon}, \tilde w_\ep,\tilde P_\ep) \in H^1_0(\Omega)^3\times H^1_0(\Omega)^3\times L^2_0(\Omega)$, we define $(\hat{u}_{\varepsilon}$, $\hat{w}_{\varepsilon},\hat P_\ep)$ by
\begin{equation}\label{uhat}
\hat{u}_{\varepsilon}(x^{\prime},y)=\tilde{u}_{\varepsilon}\left( {\varepsilon}\kappa\left(\frac{x^{\prime}}{{\varepsilon}} \right)+{\varepsilon}y^{\prime},y_3 \right),\quad
\hat{w}_{\varepsilon}(x^{\prime},y)=\tilde{w}_{\varepsilon}\left( {\varepsilon}\kappa\left(\frac{x^{\prime}}{{\varepsilon}} \right)+{\varepsilon}y^{\prime},y_3 \right)
\end{equation}
\begin{equation}\label{Phat}
\hat{P}_{\varepsilon}(x^{\prime},y)=\tilde{P}_{\varepsilon}\left( {\varepsilon}\kappa\left(\frac{x^{\prime}}{{\varepsilon}} \right)+{\varepsilon}y^{\prime},y_3 \right),
\end{equation}
a.e. $(x^{\prime},y)\in \omega\times Y$, where $\tilde u_\ep$, $\tilde w_\ep$ and $\tilde P_\ep$ are extended by zero outside $\Omega$ and the function $\kappa$ is defined by (\ref{kappa_fun}).

\begin{remark}\label{remarkCV}
For $k^{\prime}\in \mathcal{K}_{\varepsilon}$, the restrictions of $(\hat{u}_{\varepsilon}, \hat{w}_{\varepsilon}, \hat P_\ep)$  to $Y^{\prime}_{k^{\prime},a_{\varepsilon}}\times Y$ does not depend on $x^{\prime}$, whereas as a function of $y$ it is obtained from $(\tilde{u}_{\varepsilon}, \tilde w_\ep, \tilde{P}_{\varepsilon})$ by using the change of variables $y^{\prime}=\frac{x^{\prime}-a_{\varepsilon}k^{\prime}}{a_{\varepsilon}}$,
which transforms $Y_{k^{\prime},a_{\varepsilon}}$ into $Y$.
\end{remark}
Now, we get the estimates for the sequences $(\hat{u}_{\varepsilon}, \hat w_\ep, \hat{P}_{\varepsilon})$ similarly as in the proof of Lemma 4.9 in \cite{AnguianoGrau}.
 \begin{lemma}\label{estimates_hat}
 There exists a constant $C>0$ independent of $\ep$, such that $(\hat u_\ep, \hat w_\ep, \hat P_\ep)$ defined by (\ref{uhat})-(\ref{Phat}) satisfies  
 \begin{itemize}
 \item[i)] in the cases {\bf PTMP} and {\bf HTPM},
 \begin{eqnarray}
 \|\hat u_\ep\|_{L^2(\omega\times Y)^3}\leq Ca_\ep^2,& 
 \|D_{y'}\hat u_\ep\|_{L^2(\omega\times Y)^{3\times 3}}\leq Ca_\ep^2,& 
 \|\partial_{y_3}\hat u_\ep\|_{L^2(\omega\times Y)^{3}}\leq C a_\ep \ep,\label{estim_u_hat}\\
 \noame
  \|\hat w_\ep\|_{L^2(\omega\times Y)^3}\leq Ca_\ep,&
 \|D_{y'}\hat w_\ep\|_{L^2(\omega\times Y)^{3\times 3}}\leq Ca_\ep,& 
 \|\partial_{y_3}\hat w_\ep\|_{L^2(\omega\times Y)^{3}}\leq C\ep,\label{estim_w_hat}
  \end{eqnarray}
   \item[ii)] in the case {\bf VTPM},
 \begin{eqnarray}
 \|\hat u_\ep\|_{L^2(\omega\times Y)^3}\leq C\ep^2,& 
 \|D_{y'}\hat u_\ep\|_{L^2(\omega\times Y)^{3\times 3}}\leq Ca_\ep \ep,& 
 \|\partial_{y_3}\hat u_\ep\|_{L^2(\omega\times Y)^{3}}\leq C\ep^2,\label{estim_u_hat_super}\\
 \noame
  \|\hat w_\ep\|_{L^2(\omega\times Y)^3}\leq C\ep,&
 \|D_{y'}\hat w_\ep\|_{L^2(\omega\times Y)^{3\times 3}}\leq Ca_\ep,& 
 \|\partial_{y_3}\hat w_\ep\|_{L^2(\omega\times Y)^{3}}\leq C\ep,\label{estim_w_hat_super}
  \end{eqnarray}
  \end{itemize}
  and, moreover, in every cases,
  $$ \|\hat P_\ep\|_{L^2(\omega\times Y)}\leq C.\label{estim_P_hat}$$
 \end{lemma}

\paragraph{Weak variational formulation. } To finish this section, we will give the variational formulation satisfied by the functions $(\hat u_\ep,\hat w_\ep,\hat P_\ep)$, which will be useful in the following sections.

We consider $\varphi_\ep(x',y_3)=\varphi(x',x'/\ep,y_3)$ and $\psi_\ep(x',y_3)=\psi(x',x'/\ep,y_3)$ as test function in (\ref{form_var_2}) where $\varphi(x',y)$, $\psi(x',y)\in \mathcal{D}(\omega;C_{\#}^\infty(Y)^3)$, and taking into account the extension of the pressure, we have 
$$\int_{\widetilde\Omega_\ep}\nabla_{\ep}\tilde p_\ep\cdot \varphi_\ep\,dx'dy_3=\int_{\Omega}\nabla_{\ep}\tilde P_\ep\cdot \varphi_\ep\,dx'dy_3\,,$$
and the extension of $(\tilde u_\ep, \tilde w_\ep)$, we get
\begin{equation}\label{form_var_general_1}
\left\{\begin{array}{l}\displaystyle
\int_{\widetilde\Omega_\ep}D_{\ep}\tilde u_\ep:D_{\ep}\varphi_\ep\,dx'dy_3-\int_{\Omega}
\tilde P_\ep\,{\rm div}_{\ep}\varphi_\ep\,dx'dy_3
=2N^2\int_{\widetilde\Omega_\ep}{\rm rot}_{\ep}\tilde w_\ep\cdot \varphi_\ep\,dx'dy_3+\int_{\widetilde\Omega_\ep}f'\cdot \varphi_\ep'\,dx'dy_3\,,\\
\\
\displaystyle
R_M\int_{\widetilde\Omega_\ep}D_{\ep}\tilde w_\ep:D_{\ep}\psi_\ep\,dx'dy_3
+4N^2\int_{\widetilde\Omega_\ep}\tilde w_\ep\cdot \psi_\ep\,dx'dy_3=2N^2\int_{\widetilde\Omega_\ep}{\rm rot}_{\ep}\tilde u_\ep\cdot \psi_\ep\,dx'dy_3+\int_{\widetilde\Omega_\ep}g_\ep'\cdot \psi_\ep'\,dx'dy_3\,,
\end{array}\right.
\end{equation}
where $R_M$ and $g_\ep'$ depend on the case, see assumptions (\ref{estim_f_g_cases_crit_sub})-(\ref{R_M_super}).\\

Now, by the change of variables given in Remark \ref{remarkCV} (see \cite{AnguianoGrau} for more details), we obtain
\begin{equation}\label{form_var_hat_u}
\left\{\begin{array}{l}
\displaystyle{1\over a_\ep^2}\int_{\omega\times Y_f}D_{y'}\hat u_\ep':D_{y'}\varphi'\,dx'dy+{1\over \ep^2}\int_{\omega\times Y_f}\partial_{y_3}\hat u'_\ep :\partial_{y_3}\varphi'\,dx'dy\\
\noame\displaystyle
-\int_{\omega\times Y_f}\hat P_\ep{\rm div}_{x'}\varphi'\,dx'dy-{1\over a_\ep}\int_{\omega\times Y_f}\hat P_\ep{\rm div}_{y'}\varphi'\,dx'dy\\
\noame
\displaystyle
={2N^2\over a_\ep}\int_{\omega\times Y_f}{\rm rot}_{y'}\hat w_{\ep,3}\cdot \varphi'\,dx'dy+
{2N^2\over \ep}\int_{\omega\times Y_f}{\rm rot}_{y_3}\hat w_\ep'\cdot \varphi'\,dx'dy+
\int_{\omega\times Y_f}f'\cdot \varphi'\,dx'dy+O_\ep\,,
\\
\\
\\
\displaystyle
\displaystyle{1\over a_\ep^2}\int_{\omega\times Y_f}\nabla_{y'}\hat u_{\ep,3}\cdot\nabla_{y'}\varphi_3\,dx'dy+{1\over \ep^2}\int_{\omega\times Y_f}\partial_{y_3}\hat u_{\ep,3} \cdot \partial_{y_3}\varphi_3\,dx'dy
-{1\over \ep}\int_{\omega\times Y_f}\hat P_\ep\partial_{y_3}\varphi_3\,dx'dy\\
\noame
\displaystyle
={2N^2\over a_\ep}\int_{\omega\times Y_f}{\rm Rot}_{y'}\hat w_{\ep}'\, \varphi_3\,dx'dy+O_\ep\,,
\end{array}
\right.
\end{equation}
and
\begin{equation}\label{form_var_hat_w}
\left\{\begin{array}{l}
\displaystyle{R_M\over a_\ep^2}\int_{\omega\times Y_f}D_{y'}\hat w_\ep':D_{y'}\psi'\,dx'dy+{R_M\over \ep^2}\int_{\omega\times Y_f}\partial_{y_3}\hat w'_\ep :\partial_{y_3}\psi'\,dx'dy+ 4N^2\int_{\omega\times Y_f}\hat w_\ep'\cdot \psi'\,dx'dy\\
\noame
\displaystyle
={2N^2\over a_\ep}\int_{\omega\times Y_f}{\rm rot}_{y'}\hat u_{\ep,3}\cdot \psi'\,dx'dy+
{2N^2\over \ep}\int_{\omega\times Y_f}{\rm rot}_{y_3}\hat u_\ep'\cdot \psi'\,dx'dy+
\int_{\omega\times Y_f}g_\ep'\cdot \psi'\,dx'dy+O_\ep\,,
\\
\\
\\
\displaystyle
\displaystyle{R_M\over a_\ep^2} \int_{\omega\times Y_f}\nabla_{y'}\hat w_{\ep,3}\cdot \nabla_{y'}\psi_3\,dx'dy+{R_M\over \ep^2}\int_{\omega\times Y_f}\partial_{y_3}\hat w_{\ep,3} :\partial_{y_3}\psi_3\,dx'dy
+4N^2\int_{\omega\times Y_f}\hat w_{\ep,3}\cdot \psi_3\,dx'dy\\
\noame
\displaystyle={2N^2\over a_\ep}\int_{\omega\times Y_f}{\rm Rot}_{y'}\hat u_{\ep}'\, \psi_3\,dx'dy+O_\ep\,.
\end{array}
\right.
\end{equation}

When $\ep$ tends to zero, we obtain for $(\hat u_\ep, \hat w_\ep, \hat P_\ep)$ different asymptotic behaviors depending on the cases {\bf PTPM}, {\bf HTPM} and {\bf VTPM}. We will analyze them in the next sections.

\section{Proportionally Thin Porous Medium (PTPM)}\label{PTPM}
It corresponds to the critical case when the cylinder height is proportional to the interspatial distance, with $\lambda$ the proportionality constant, that is $a_\ep\approx\ep$, with $a_\ep/\ep\to \lambda$, $0<\lambda<+\infty$. 
%

Let us introduce some notation which will be useful along this section. For a vectorial function $v=(v',v_3)$ and a scalar function $w$, we introduce the operators $D_\lambda$, $\nabla_\lambda$, ${\rm div}_\lambda$ and ${\rm rot}_\lambda$ by 
\begin{eqnarray}
&(D_{\lambda}v)_{ij}=\partial_{x_j}v_i\hbox{ for }i=1,2,3,\ j=1,2,\quad (D_{\lambda}v)_{i,3}=\lambda \partial_{y_3}v_i\hbox{ for }i=1,2,3,\nonumber&\\
\noame
&\Delta_\lambda v=\Delta_{y'}v+\lambda^2\partial_{y_3}^2 v,\quad \nabla_{\lambda} w=( \nabla_{y'}w,\lambda \partial_{y_3} w)^t,&\nonumber\\
\noame
&{\rm div}_\lambda v={\rm div}_{y'}v'+\lambda\partial_{y_3} v_3,\quad {\rm rot}_\lambda v=({\rm rot}_{y'}v_3+
\lambda{\rm rot}_{y_3} v',{\rm Rot}_{y'} v')\,,&\nonumber
\end{eqnarray}
where ${\rm rot}_{y'}$, ${\rm rot}_{y_3}$ and ${\rm Rot}_{y'}$ are defined in (\ref{def_rot}).
Next, we give some compactness results about the behavior of the extended sequences $(\tilde u_\ep, \tilde w_\ep, \tilde P_\ep)$ and the unfolding functions $(\hat u_\ep, \hat w_\ep, \hat P_\ep)$ satisfying the {\it a priori} estimates given in Lemmas \ref{lemma_estimates},  \ref{lemma_est_P} and \ref{estimates_hat} respectively.
\begin{lemma}\label{lem_asymp_crit}
For a subsequence of $\ep$ still denote by $\ep$, we have that 
\begin{itemize}
\item[i)] {\it (Velocity)} there exist $\tilde u\in H^1_0(0,1;L^2(\omega)^3)$ with $\tilde u_3=0$  and  $\hat u\in L^2(\omega;H^1_{0,\#}(Y))^3$ with $\hat u=0$ on $\omega\times Y_s$, such that $\int_{Y}\hat u(x',y)dy=\int_0^1\tilde u(x',y_3)\,dy_3$ with $\int_{Y}\hat u_3\,dy=0$ and moreover 
\begin{eqnarray}
&\displaystyle a_\ep^{-2}\tilde u_\ep\rightharpoonup (\tilde u',0)\hbox{  in  }H^1(0,1;L^2(\omega)^3),\quad 
a_\ep^{-2}\hat u_\ep\rightharpoonup \hat u\hbox{  in  }L^2(\omega;H^1(Y)^3),&\label{conv_u_crit}\\
\noame
&\displaystyle{\rm div}_{x'}\left(\int_0^{1}\tilde u'(x',y_3)\,dy_3\right)=0\hbox{  in  }\omega,\quad 
\left(\int_0^1\tilde u'(x',y_3)\,dy_3\right)\cdot n=0\hbox{  in  }\partial\omega\,,&\label{div_x_crit}\\
\noame
&\displaystyle{\rm div}_\lambda\hat u=0\hbox{  in  }\omega\times Y,\quad 
{\rm div}_{x'}\left(\int_{Y}\hat u'(x',y)\,dy\right)=0\hbox{  in  }\omega,\quad \left(\int_{Y}\hat u'(x',y)\,dy\right)\cdot n=0
\hbox{  in  }\partial\omega\,,&\label{div_y_crit}
\end{eqnarray}
\item[ii)] {\it (Microrotation)} there exist $\tilde w\in H^1_0(0,1;L^2(\omega)^3)$ with  $\tilde w_3=0$ and  $\hat w\in L^2(\omega;H^1_{0,\#}(Y))^3$ with $\hat w=0$ on $\omega\times Y_s$, such that $\int_{Y}\hat w(x',y)dy=\int_0^{1}\tilde w(x',y_3)\,dy_3$ with $\int_{Y}\hat w_3\,dy=0$ and moreover 
\begin{eqnarray}
&\displaystyle a_\ep^{-1}\tilde w_\ep\rightharpoonup (\tilde w',0)\hbox{  in  }H^1(0,1;L^2(\omega)^3),\quad 
a_\ep^{-1}\hat w_\ep\rightharpoonup \hat w\hbox{  in  }L^2(\omega;H^1(Y)^3),&\label{conv_w_crit}
\end{eqnarray}
\item[iii)] {\it (Pressure)} there exists a function $\tilde P\in L^2_0(\Omega)$, independent of $y_3$, such that
\begin{eqnarray}
&\displaystyle\tilde P_\ep\to \tilde P\hbox{  in  }L^2(\Omega),\quad \hat P_\ep\to \tilde P\hbox{  in  }L^2(\omega\times Y).& \label{conv_P_crit}
\end{eqnarray}
\end{itemize}
\end{lemma}
{\bf Proof. }  The proof of this result for the velocity is obtained by arguing similarly to Section 5 in \cite{AnguianoGrau}.

The proof of the results for the microrotation is analogous to the ones of the velocity, except to prove that $\tilde w_3=0$. To do this, we consider as test function $\psi_\ep(x',y_3)=(0,0,a_\ep^{-1}\psi_3)$ in the variational formulation (\ref{form_var_general_1}), and we get
$$
\begin{array}{l}
\displaystyle a_\ep R_c\int_{\Omega}\nabla_{x'}\tilde w_{\ep,3}\cdot \nabla_{x'}\psi_3\,dx'dy_3+a_\ep\ep^{-2}R_c\int_\Omega
\partial_{y_3}\tilde w_{\ep,3}\,\partial_{y_3}\psi_3\,dx'dy_3+4N^2 a_\ep^{-1}\int_\Omega\tilde w_{\ep,3}\psi_3\,dx'dy_3\\
\noame
\displaystyle=2N^2 a_\ep^{-1}\int_{\Omega}{\rm Rot}_{x'}\tilde u'_\ep\,\psi_3\,dx'dy_3\,.
\end{array}
$$
Passing to the limit by using convergences of $\tilde u_\ep$ and $\tilde w_\ep$ given in (\ref{conv_u_crit}) and (\ref{conv_w_crit}), we get
$$\lambda^2R_c\int_\Omega \partial_{y_3}\tilde w_3\,\partial_{y_3}\psi_3\,dx'dy_3+4N^2\int_\Omega \tilde w_3\,\psi_3\,dx'dy_3=0\,,$$
and taking into account that $\tilde w_3=0$ on $y_3=\{0,1\}$,  it is easily deduced that $\tilde w_3=0$ a.e. in $\Omega$.

We finish with the proof for the pressure. Estimate (\ref{estim_P_hat}) implies, up to a subsequence, the existence of $\tilde P\in L^2_0(\Omega)$ such that 
\begin{equation}\label{conv_p_weak}
\tilde P_\ep\rightharpoonup \tilde P\quad \hbox{  in  }L^2(\Omega).
\end{equation} Also, from $\|\nabla_{\ep}\tilde P_\ep\|_{L^2(\Omega)^3}\leq C$, by noting that $\partial_{y_3}\tilde P_\ep/\ep$ also converges weakly in $H^{-1}(\Omega)$, we obtain $\partial_{y_3}\tilde P=0$ and so $\tilde P$ is independent of $y_3$.\\

Next, following  \cite{Tartar}, we prove that the convergence of the pressure is in fact strong. As $\tilde u_3=0$, $\partial_{y_3}\tilde P_\ep$ tends to zero and $\tilde p$ only depends on $x'$, let $\sigma_\ep(x',y_3)=(\sigma_\ep'(x'),0)\in H^1_0(\omega)^3$ be such that 
\begin{equation}\label{strong_p_1}
\sigma_\ep\rightharpoonup \sigma\quad\hbox{in }H^1_0(\omega)^3. 
\end{equation}
Then, we have
$$\left|<\nabla_{x'}\tilde P_\ep,\sigma_\ep>_{\Omega}-<\nabla_{x'}\tilde P,\sigma>_{\Omega}\right|\leq
\left|<\nabla_{x'}\tilde P_\ep,\sigma_\ep-\sigma>_{\Omega}\right|+\left|<\nabla_{x'}\tilde P_\ep-\nabla_{x'}\tilde P,\sigma>_{\Omega}\right|.$$
On the one hand, using convergence (\ref{conv_p_weak}), we have 
$$\left|<\nabla_{x'} \tilde P_\ep-\nabla_{x'}\tilde P,\sigma>_\Omega\right|=\left|\int_\Omega\left(\tilde P_\ep-\tilde P\right)\,{\rm div}_{x'}\sigma'\,dx\right|\to 0,\quad \hbox{as }\ep\to 0\,.$$
On the other hand, from (\ref{extension_1}) and following the proof of Lemma \ref{lemma_est_P}, we have that
$$\begin{array}{ll}
\left|<\nabla_{x'}\tilde P_\ep,\sigma_\ep-\sigma>_\Omega\right|=& \left|<\nabla_{x'}\tilde P_\ep,\tilde R^\ep(\sigma_\ep'-\sigma')>_{\widetilde\Omega_\ep}\right|\\
\noame
& \leq C\left(\|\sigma'_\ep-\sigma'\|_{L^2(\omega)^3}+a_\ep\|D_{x'}(\sigma'_\ep-\sigma')\|_{L^2(\omega)^3}\right)\to 0\quad\hbox{as }\ep\to 0,
\end{array}$$
by virtue of (\ref{strong_p_1}) and the Rellich theorem. This implies that $\nabla_{x'}\tilde P_\ep\to \nabla_{x'}\tilde P$ strongly in $H^{-1}(\Omega)^3$, which together the Ne${\breve{\rm c}}$as inequality, implies the strong convergence of the pressure $\tilde P_\ep$ given in (\ref{conv_P_crit}).  We finish the proof by proving that $\hat P_\ep$ also converges strongly to $\tilde P$. The existence of $\hat q\in L^2(\omega\times Y)$ such that   $\hat P_\ep$ converges weakly to $\hat q$  is a consequence of estimate (\ref{estim_P_hat}). It remains only to prove that $\hat q$ does not depend on the microscopic variable $y$. To do this, we choose as test function $\varphi_\ep(x',y)=(\ep\varphi',\eta_\ep\varphi_3)$ in (\ref{form_var_hat_u}). Taking into account the estimates (\ref{estim_u_hat})-(\ref{estim_P_hat}) and passing to the limit when $\ep$ tends to zero by using convergence (\ref{conv_P_crit}), we have
$$\int_{\omega\times Y}\hat P{\rm div}_{y}\varphi\,dx'dy=0\,,$$
which shows that $\hat P$ does not depend on $y$, which ends the proof.
\par\hfill$\square$\\


Unsing previous convergences, in the following theorem we give the homogenized system satisfied by $(\hat u, \hat w, \tilde P)$.
\begin{theorem}
In the case {\bf PTPM},   the sequence $(a_\ep^{-2}\hat u_\ep, a_\ep^{-1}\hat w_\ep)$ converges weakly to $(\hat u,\hat w)$ in $L^2(\omega;H^1(Y)^3)\times L^2(\omega;H^1(Y)^3)$ and  $\hat P_\ep$ converges strongly to $\tilde P$ in $L^2(\omega)$, where $(\hat u,\hat w, \tilde P)\in L^2(\omega;H^1_{0,\#}(Y)^3)\times L^2(\omega;H^1_{0,\#}(Y)^3)\times (L^2_0(\omega)\cap H^1(\omega))$, with $\int_Y\hat u_3\,dy=\int_Y\hat w_3\,dy=0$, is the unique solution of the following homogenized system
\begin{equation}\label{hom_system_crit}
\left\{\begin{array}{rl}
\displaystyle
-\Delta_\lambda \hat u+\nabla_\lambda \hat q=2N^2 {\rm rot}_\lambda\hat w+f'(x')-\nabla_{x'}\tilde P(x')&\hbox{ in }\omega\times Y_f,\\
\noame
{\rm div}_\lambda\hat u=0&\hbox{ in }\omega\times Y_f,\\
\noame
-R_c\Delta_\lambda \hat w+4N^2  \hat w=2N^2 \, {\rm rot}_\lambda\hat u+  g'(x')&\hbox{ in }\omega\times Y_f,\\
\noame
\hat u=\hat w=0&\hbox{ in }\omega\times Y_s,\\
\noame
\displaystyle{\rm div}_{x'}\left(\int_{Y}\hat u'(x',y)\,dy\right)=0&\hbox{ in }\omega,\\
\noame
\displaystyle\left(\int_{Y}\hat u'(x',y)\,dy\right)\cdot n=0&\hbox{ on }\partial\omega,\\
\noame
\hat u(x',y), \hat w(x',y), \hat q(x',y)\quad Y'-\hbox{periodic}.&
\end{array}\right.
\end{equation}
\end{theorem}
{\bf Proof. } For every $\varphi\in \mathcal{D}(\omega;C_{\#}^\infty(Y)^3)$ with ${\rm div}_\lambda\varphi=0$ in $\omega\times Y$ and ${\rm div}_{x'}(\int_Y \varphi'\,dy)=0$ in $\omega$, we choose $\varphi_\ep=(\varphi',\lambda(\ep/a_\ep)\varphi_3)$ in (\ref{form_var_hat_u}). Taking into account that thanks to ${\rm div}_\lambda\varphi=0$ in $\omega\times Y$, we have that 
$${1\over a_\ep}\int_{\omega\times Y} \hat P_\ep ({\rm div}_{y'}\varphi'+\lambda\partial_{y_3}\varphi_3)\,dx'dy=0\,.$$
Thus, passing to the limit using the convergences (\ref{conv_u_crit}), (\ref{conv_w_crit}), (\ref{conv_P_crit}) and $\lambda(\ep/\eta_\ep)\to 1$, we obtain
\begin{equation}\label{form_varl_limit_crit_con_p}
\begin{array}{l}
\displaystyle\int_{\omega\times Y_f}D_\lambda\hat u:D_\lambda  \varphi\,dx'dy-\int_{\omega\times Y}\tilde P\,{\rm div}_{x'}\varphi'\,dx'dy\\
\noame
\displaystyle
=2N^2 \int_{\omega\times Y_f}\left({\rm rot}_{y'}\hat w_3\cdot \varphi'+\lambda{\rm rot}_{y_3}\hat w'\cdot \varphi'+ {\rm Rot}_{y'}\hat w'\,\varphi_3\right)dx'dy+\int_{\omega\times Y_f}f'\cdot \varphi'\,dx'dy\,.
\end{array}
\end{equation}
Since $\tilde P$ does not depend on $y$ and ${\rm div}_{x'}\int_{Y}\varphi'\,dy=0$ in $\omega$, we have that the second term is zero, and 
so we get
\begin{equation}\label{form_var_limit_crit}\int_{\omega\times Y_f}D_\lambda\hat u:D_\lambda\varphi\,dx'dy=2N^2\int_{\omega\times Y_f}{\rm rot}_{\lambda}\hat w\cdot \varphi\,dx'dy+\int_{\omega\times Y_f}f'\cdot \varphi'\,dx'dy\,.
\end{equation}

Next, for every $\psi\in \mathcal{D}(\omega;C_\#^\infty(Y)^3)$, we choose $\psi_\ep=a_\ep^{-1}\psi$ in (\ref{form_var_hat_w}) with $g_\ep$ and $R_M$ satisfying (\ref{estim_f_g_cases_crit_sub}) and (\ref{R_M}). Then, passing to the limit using convergences (\ref{conv_u_crit}) and (\ref{conv_w_crit}), we get
\begin{equation}\label{form_var_limit_crit_w}
\begin{array}{l}\displaystyle R_c\int_{\omega\times Y_f}D_\lambda\hat w:D_\lambda\psi\,dx'dy+4N^2 \int_{\omega\times Y_f}\hat w\cdot \psi\,dx'dy =2N^2  \int_{\omega\times Y_f}{\rm rot}_{\lambda}\hat u\cdot \psi\,dx'dy+  \int_{\omega\times Y_f}g'\cdot \psi'\,dx'dy\,.
\end{array}
\end{equation}
By density (\ref{form_var_limit_crit}) holds for every function $\varphi$ in the Hilbert space $V$ defined by
$$
V=\left\{\begin{array}{l}
\varphi(x',y)\in L^2(\omega;H^1_{0,\#}(Y)^3) \hbox{ such that }\\
\noame
\displaystyle {\rm div}_{x'}\left(\int_{Y_f}\varphi(x',y)\,dy \right)=0\hbox{ in }\omega,\quad \left(\int_{Y_f}\varphi(x',y)\,dy \right)\cdot n=0\hbox{ on }\partial\omega\\\noame
{\rm div}_\lambda\varphi(x',y)=0\hbox{ in }\omega\times Y_f,\quad \varphi(x',y)=0\hbox{ in }\omega\times Y_s
\end{array}\right\}\,,
$$
and (\ref{form_var_limit_crit_w}) in $W=\{\psi(x',y)\in L^2(\omega;H^1_{0,\#}(Y)^3)\,:\, \psi(x',y)=0\hbox{ in }\omega\times Y_s\}$.

From Theorem 2.4.2 in \cite{Luka}, the variational formulation (\ref{form_var_limit_crit})-(\ref{form_var_limit_crit_w}) admits a unique solution $(\hat u, \hat w)$ in $V\times W$.

From Lemma 2.4.1 in \cite{Luka} (see also \cite{Allaire0}), the orthogonal of $V$ with respect to the usual scalar product in $L^2(\omega\times Y)$ is made of gradients of the form $\nabla_{x'}q(x')+\nabla_{\lambda}\hat q(x',y)$, with $q(x')\in L^2_0(\omega)$ and $\hat q(x',y)\in L^2(\omega;H^1_\#(Y))$. Therefore, by integration by parts, the variational formulations (\ref{form_var_limit_crit})-(\ref{form_var_limit_crit_w})  are equivalent to the homogenized system (\ref{hom_system_crit}). It remains to prove that the pressure $\tilde P(x')$, arising as a Lagrange multiplier of the incompressibility constraint ${\rm div}_{x'}(\int_Y \hat u(x',y)dy)=0$, is the same as the limit of the pressure $\hat P_\ep$. This can be easily done by considering in equation (\ref{form_var_hat_u}) a test function with ${\rm div}_\lambda$ equal to zero, and  obtain  the variational formulation (\ref{form_varl_limit_crit_con_p}). Since $2N^2{\rm rot}_\lambda\hat w+ f'\in L^2(\omega\times Y)^3$ we deduce that $\tilde P\in H^1(\omega)$.

Finally, since from Lemma 2.4.1 in \cite{Luka}  we have that  (\ref{hom_system_crit}) admits a unique solution, and then the complete sequence $(a_\ep^{-2}\hat u_\ep, a_\ep^{-1}\hat w_\ep, \hat P_\ep)$ converges to the solution $(\hat u(x',y), \hat w(x',y), \tilde P(x'))$. 
\par\hfill$\square$

Let us define the local problems which are useful to eliminate the variable $y$ of the previous homogenized problem and then obtain a Darcy equation for the pressure $\tilde P$.

For every $i,k= 1, 2$ and $0< \lambda<+\infty$, we consider the following local micropolar problems 
\begin{equation}\label{local_problems_crit}
\left\{\begin{array}{rl}\displaystyle
-\Delta_\lambda u^{i,k}+\nabla_{\lambda}\pi^{i,k}-2N^2{\rm rot}_\lambda w^{i,k}=e_i\delta_{1k}&\hbox{ in }Y_f,\\
\noame
{\rm div}_\lambda u^{i,k}=0&\hbox{ in }Y_f,\\
\noame
-R_c\Delta_\lambda w^{i,k}+4N^2  w^{i,k}-2N^2 \,{\rm rot}_\lambda u^{i,k}= \,e_i\delta_{2k}&\hbox{ in }Y_f,\\
\noame
u^{i,k}=w^{i,k}=0&\hbox{ in }\omega\times Y_s,\\
\noame
\displaystyle\int_{Y_f}u_3^{i,k}(y)dy=\int_{Y_f}w_3^{i,k}(y)dy=0&\hbox{ in }\omega,\\
\noame
u^{i,k}(y), w^{i,k}(y),\pi^{i,k}(y)\quad Y'-\hbox{periodic}.
\end{array}\right.
\end{equation}
It is known (see Lemma 2.5.1 in \cite{Luka}) that there exist a unique solution $(u^{i,k},w^{i,k},\pi^{i,k})\in H^1_{0,\#}(Y_f)^3\times H^1_{0,\#}(Y_f)^3\times L^2_0(Y_f)$ of problem (\ref{local_problems_crit}), and moreover $\pi^{i,k}\in H^1(Y_f)$.

Now, we give the main result concerning the homogenized flow.  

\begin{theorem}\label{them_main_crit}
Let $(\hat u,\hat w,\tilde P)\in L^2(\omega;H^1_{0,\#}(Y)^3)\times L^2(\omega;H^1_{0,\#}(Y)^3)\times (L_0^2(\omega)\cap H^1(\omega))$ be the unique weak solution of problem (\ref{hom_system_crit}). Then, the extensions $(a_\ep^{-2}\tilde u_\ep,a_\ep^{-1}\tilde w_\ep)$ and  $\tilde P_\ep$ of the solution of problem (\ref{system_2})-(\ref{bc_system_2}) converge weakly to $(\tilde u,\tilde w)$   in $H^1(0, 1;L^2(\omega)^3)\times H^1(0,1;L^2(\omega)^3)$ and strongly to $\tilde P$ in $L^2(\omega)$ respectively, with $\tilde u_3=\tilde w_3=0$. Moreover, defining $\widetilde U(x')=\int_0^{1}\tilde u(x',y_3)\,dy_3$ and $\widetilde W(x')=\int_0^{1}\tilde w(x',y_3)\,dy_3$, it holds 
\begin{equation}\label{Darcy_law_u_w_crit}
\begin{array}{ll}
\widetilde U'(x')=K^{(1)}_\lambda\left(f'(x')-\nabla_{x'}\tilde P(x')\right)+K^{(2)}_\lambda g(x'),&\quad \widetilde U_3(x')=0\quad \hbox{ in }\omega,\\
\noame
\widetilde W'(x')=L^{(1)}_\lambda\left(f'(x')-\nabla_{x'}\tilde P(x')\right)+L^{(2)}_\lambda g(x'),&\quad \widetilde W_3(x')=0\quad \hbox{ in }\omega,
\end{array}
\end{equation}
where $K^{(k)}_\lambda$, $L^{(k)}_\lambda\in \mathbb{R}^{2\times 2}$, $k=1,2$, are matrices with coefficients
$$\left(K^{(k)}_\lambda\right)_{ij}=\int_{Y_f} u^{i,k}_j(y)\,dy,\quad \left(L^{(k)}_\lambda\right)_{ij}=\int_{Y_f} w^{i,k}_j(y)\,dy,\quad i,j=1,2,$$
where $u^{i,k}$, $w^{i,k}$ are the solutions of the local micropolar problems defined in (\ref{local_problems_crit}).

Here, $\tilde P\in H^1(\omega)\cap L^2_0(\omega)$ is the unique solution of the 2D Darcy equation
\begin{equation}\label{Darcy_law_P_crit}
\left\{\begin{array}{l}
{\rm div}_{x'}\left( K^{(1)}_\lambda\left(f'(x')-\nabla_{x'}\tilde P(x')\right)+K^{(2)}_\lambda g(x')\right)=0\quad \hbox{ in }\omega,\\
\noame
\left( K^{(1)}_\lambda\left(f'(x')-\nabla_{x'}\tilde P(x')\right)+K^{(2)}_\lambda g(x')\right)\cdot n=0\quad \hbox{ in }\partial\omega.
\end{array}\right.
\end{equation}
\end{theorem}
{\bf Proof. }We eliminate the microscopic variable $y$ in the effective problem (\ref{hom_system_crit}). To do that, we consider the following identification
\begin{eqnarray}
\hat u(x',y)=\sum_{i=1}^2\left[\left(f_i(x')-\partial_{x_i}\tilde P(x')\right) u^{i,1}(y)+ g_i(x') u^{i,2}(y)\right],\nonumber\\
\noame\displaystyle \hat w(x',y)=\sum_{i=1}^2\left[\left(f_i(x')-\partial_{x_i}\tilde P(x')\right) w^{i,1}(y)+ g_i(x') w^{i,2}(y)\right],\label{idendity_for_y}\\
\noame
\hat q(x',y)=\sum_{i=1}^2\left[\left(f_i(x')-\partial_{x_i}\tilde P(x')\right) \pi^{i,1}(y) + g_i(x') \pi^{i,2}(y)\right],\nonumber
\end{eqnarray}
and thanks to the identity $\int_{Y_f} \hat \varphi(x',y)\,dy=\int_0^{1}\tilde \varphi(x',y_3)\,dy_3$ with $\int_{Y_f}\hat \varphi_3\,dy=0$ satisfied by  velocity and microrotation  given in Lemma \ref{lem_asymp_crit}, we deduce that $\widetilde U$ and $\widetilde W$ are given by (\ref{Darcy_law_u_w_crit}).

Finally, the divergence condition with respect to the variable $x'$ given in  (\ref{hom_system_crit}) together with the expression of $\widetilde U'(x')$ gives (\ref{Darcy_law_P_crit}),  which has a unique solution since $K^{(1)}_\lambda$ is positive definite and then the whole sequence converges, see Part III - Theorem 2.5.2  in \cite{Luka}.
\par\hfill$\square$

\begin{remark} We observe that when $N$ is identically zero, taking into account the linear momentum equations from (\ref{hom_system_crit}), we can deduce that the Darcy equation (\ref{Darcy_law_P_crit})  agrees with the ones obtained in \cite{AnguianoGrau2, Fabricius} in the case {\bf PTPM}. 
\end{remark}

\section{The homogeneously thin porous medium (HTPM)}\label{HTPM}
It corresponds to the case when the cylinder height is much larger than interspatial distance, i.e. $a_\ep\ll \ep$ which is equivalent to $\lambda=0$.

%
Next, we give some compactness results about the behavior of the extended sequences $(\tilde u_\ep,\tilde w_\ep,\tilde P_\ep)$ and the unfolding functions $(\hat u_\ep,\hat w_\ep, \hat P_\ep)$ by using the {\it a priori} estimates given in Lemmas \ref{lemma_estimates} and \ref{lemma_est_P}, and Lemma \ref{estimates_hat}, respectively.

\begin{lemma}\label{lem_asymp_HTPM}
For a subsequence of $\ep$ still denoted by $\ep$, there exist the following functions:
\begin{itemize}
\item[i)] (Velocity)  there exist $\tilde u\in L^2(\Omega)^3$, with $\tilde u_3=0$ and  $\hat u\in   L^2(\Omega; H^1_\#(Y')^3)$ with $\hat u=0$ on $\omega\times Y_s$, such that $\int_{Y_f}\hat u(x',y)dy=\int_0^1\tilde u(x',y_3)dy_3$ with $\int_{Y_f}\hat u_3(x',y)dy=0$, $\hat u_3$ independent of $y_3$ and moreover
\begin{eqnarray}
&\displaystyle a_\ep^{-2}\tilde u_\ep\rightharpoonup (\tilde u',0)\hbox{  in  } H^1(0,1;L^2(\omega)^3),\quad 
a_\ep^{-2}\hat u_\ep\rightharpoonup \hat u\hbox{  in  } L^2(\Omega;H^1(Y')^3),&\label{conv_u_HTPM}\\
\noame
&\displaystyle{\rm div}_{x'}\left(\int_0^1 \tilde u'(x',y_3)\,dy_3\right) =0\hbox{  in  }\omega,\quad 
\left(\int_0^1 \tilde u'(x',y_3)\,dy_3\right) \cdot n=0\hbox{  in  }\partial\omega\,,&\label{div_x_HTPM}\\
\noame
&\displaystyle{\rm div}_{y'}\hat u'=0\hbox{  in  }\omega\times Y_f,\quad 
{\rm div}_{x'}\left(\int_{Y_f}\hat u'(x',y)\,dy\right)=0\hbox{  in  }\omega,\quad \left(\int_{Y_f}\hat u(x',y)\,dy\right)\cdot n=0
\hbox{  in  }\partial\omega\,,&\label{div_y_HTPM}
\end{eqnarray}

\item[ii)] (Microrotation) there exist  $\tilde w\in  L^2(\Omega)^3$ with $\tilde w_3=0$ and $\hat w\in L^2(\Omega; H^1_\#(Y')^3)$ with  $\hat w=0$ on $\omega\times Y_s$, such that $\int_{Y_f}\hat w(x',y)dy=\int_0^{1}\tilde w(x',y_3)\,dy_3$ with  $\int_{Y_f}\hat w_3(x',y)\,dy=0$,  $\hat u_3$ independent of $y_3$  and moreover 
\begin{eqnarray}
&\displaystyle a_\ep^{-1}\tilde w_\ep\rightharpoonup (\tilde w',0)\hbox{  in  }H^1(0,1;L^2(\omega)^3),\quad 
 a_\ep^{-1}\hat w_\ep\rightharpoonup \hat w\hbox{  in  }L^2(\Omega; H^1(Y')^3),&\label{conv_w_HTPM}
\end{eqnarray}
\item[iii)] (Pressure) there exists $\tilde P\in L^2_0(\Omega)$ independent of $y_3$, such that 
\begin{eqnarray}
&\displaystyle\tilde P_\ep\to \tilde P\hbox{  in  }L^2(\omega),\quad \hat P_\ep\to \tilde P\hbox{  in  }L^2(\omega).& \label{conv_P_HTPM}
\end{eqnarray}
\end{itemize}
\end{lemma}

{\bf Proof. } The proof of this result is obtained by arguing similarly to Section 5 in \cite{AnguianoGrau} and Lemma \ref{lem_asymp_crit} of the present paper.\par\hfill$\square$

Using previous convergences, in the following theorem we give the homogenized system satisfied by $(\hat u, \hat w, \tilde P)$.
\begin{theorem}
In the case {\bf HTPM},  the sequence $(a_\ep^{-2}\hat u_\ep, a_\ep^{-1}\hat w_\ep)$ converges weakly to $(\hat u,\hat w)$ in $L^2(\Omega; H^1(Y')^3)\times L^2(\Omega; H^1(Y')^3)$ and  $\hat P_\ep$ converges strongly to $\tilde P$ in $L^2(\omega)$, where $(\hat u,\hat w, \tilde P)\in L^2(\Omega; H^1_\#(Y')^3)\times L^2(\Omega; H^1_\#(Y')^3)\times (L^2_0(\omega)\cap H^1(\omega))$ with $\hat u_3$ and $\hat w_3$ independent of $y_3$ and   $\int_{Y_f'}\hat u_3(x',y')\,dy'=\int_{Y_f'}\hat w_3(x',y')\,dy'=0$. Moreover, defining  $\hat U=\int_0^1\hat u(x',y)dy_3$,  $\hat W=\int_0^1\hat u(x',y)dy_3$, we have that $(\hat U, \hat W)\in L^2(\omega; H^1_\#(Y')^3)\times L^2(\omega; H^1_\#(Y')^3)$ is the unique solution of the following homogenized system
\begin{equation}\label{hom_system_HTPM}
\left\{\begin{array}{rl}
\displaystyle
-\Delta_{y'} \hat W'+\nabla_{y'} \hat q=2N^2 {\rm rot}_{y'}\hat W_3+f'(x')-\nabla_{x'}\tilde P(x')&\hbox{ in }\omega\times Y'_f,\\
\noame
-\Delta_{y'} \hat W_3=2N^2 {\rm Rot}_{y'}\hat W'&\hbox{ in }\omega\times Y'_f,\\
\noame
-R_c\Delta_{y'} \hat W'+4N^2 \hat W'=2N^2 {\rm rot}_{y'}\hat U_3+g'(x')&\hbox{ in }\omega\times Y'_f,\\
\noame
-R_c\Delta_{y'} \hat W_3+4N^2 \hat W_3=2N^2 {\rm Rot}_{y'}\hat U'&\hbox{ in }\omega\times Y'_f,\\
\noame
{\rm div}_{y'}\hat U'=0&\hbox{ in }\omega\times Y'_f,\\
\noame
\hat U'=\hat W'=0&\hbox{ in }\omega\times Y'_s,\\
\noame
\displaystyle{\rm div}_{x'}\left(\int_{Y'_f}\hat U'(x',y')\,dy'\right)=0&\hbox{ in }\omega,\\
\noame
\displaystyle\left(\int_{Y'_f}\hat U'(x',y')\,dy'\right)\cdot n=0&\hbox{ on }\partial\omega,\\
\noame
\hat U(x',y'),\hat w(x',y'), \hat q(x',y')\quad Y'-\hbox{periodic}.&
\end{array}\right.
\end{equation}
\end{theorem}

{\bf Proof. } We choose $\varphi\in \mathcal{D}(\omega;C_{\#}^\infty(Y)^3)$ with ${\rm div}_{y'}\varphi'=0$ in $\omega\times Y$, ${\rm div}_{x'}(\int_Y \varphi'\,dy)=0$ in $\omega$ and  $\varphi_3$ independent of $y_3$ in (\ref{form_var_hat_u}). Taking into account that thanks to ${\rm div}_{y'}\varphi'=0$ in $\omega\times Y_f$ and $\varphi_3$ independent of $y_3$, we have that 
$${1\over a_\ep}\int_{\omega\times Y} \hat P_\ep {\rm div}_{y'}\varphi'\,dx'dy=0\quad\hbox{and}\quad {1\over \ep}\int_{\omega\times Y} \partial_{y_3}\hat P_\ep \partial_{y_3}\varphi_3\,dx'dy=0\,.$$
Thus, passing to the limit using the convergences (\ref{conv_u_HTPM}), (\ref{conv_w_HTPM}), (\ref{conv_P_HTPM}), $a_\ep/\ep\to 0$ and using in the limit that $\tilde P$ does not depend on $y$ and ${\rm div}_{x'}(\int_Y \varphi'\,dy)=0$, we obtain
\begin{equation}\label{form_varl_limit_HTPM_u}
\left\{\begin{array}{l}
\displaystyle\int_{\omega\times Y_f}D_{y'}\hat u':D_{y'} \varphi'\,dx'dy=2N^2 \int_{\omega\times Y_f}{\rm rot}_{y'}\hat w_3\cdot \varphi'\,dx'dy+\int_{\omega\times Y_f}f'\cdot \varphi'\,dx'dy\\
\noame
\displaystyle\int_{\omega\times Y_f}\nabla_{y'}\hat u_3:\nabla_{y'} \varphi_3\,dx'dy=2N^2 \int_{\omega\times Y_f}{\rm Rot}_{y'}\hat w'\,\varphi_3\,dx'dy\,.
\end{array}\right.
\end{equation}

Next, for every $\psi\in \mathcal{D}(\omega;C_\#^\infty(Y)^3)$ with $\psi_3$ independent of $y_3$, we choose $\psi_\ep=a_\ep^{-1}\psi$ in (\ref{form_var_hat_w}) taking into account that $g_\ep$ and $R_M$ satisfy  (\ref{estim_f_g_cases_crit_sub}) and (\ref{R_M}). Then, passing to the limit using convergences (\ref{conv_u_HTPM}) and (\ref{conv_w_HTPM}), we get
\begin{equation}\label{form_var_limit_HTPM_w}
\left\{\begin{array}{l}\displaystyle R_c\int_{\omega\times Y_f}D_{y'}\hat w':D_{y'}\psi'\,dx'dy+4N^2\int_{\omega\times Y_f}\hat w'\cdot \psi'\,dx'dy={2N^2}\int_{\omega\times Y_f}{\rm rot}_{y'}\hat u_3\cdot \psi'\,dx'dy+\int_{\omega\times Y_f}g'(x')\cdot \psi'\,dx'dy\,,\\
\noame
\displaystyle R_c\int_{\omega\times Y_f}\nabla_{y'}\hat w_3:\nabla_{y'}\psi_3\,dx'dy+4N^2\int_{\omega\times Y_f}\hat w_3\cdot \psi_3\,dx'dy=2N^2\int_{\omega\times Y_f}{\rm Rot}_{y'}\hat u'\cdot \psi_3\,dx'dy.
\end{array}\right.
\end{equation}
We take into account that there is no $y_3$-dependence in the obtained  variational formulation. For that, we can consider $\varphi$,$\psi$ independent of $y_3$, which implies that $(\hat U, \hat W)$ satisfies the same variational formulation with integral in $\omega\times Y'_f$. By density, we can deduce that the variational formulation for  $(\hat U, \hat W)$ is equivalent to problem (\ref{hom_system_HTPM}).
\par\hfill$\square$

The local problems to eliminate the variable $y$ of the previous homogenized problem can be defined by using the local system (\ref{local_problems_crit}) with $\lambda=0$. In that case, since there is no $y_3$ dependence, then it is correct to consider all equations in 2D domain $Y_f'$ instead of $Y_f$. Therefore, for every $i,k= 1, 2$, we consider $(u^{i,k},w^{i,k},\pi^{i,k})\in H^1(Y'_f)^3\times H^1(Y'_f)^3\times (H^1(Y_f')\cap L^2_0(Y_f'))$ the unique solutions of  the following local micropolar problems 
\begin{equation}\label{local_problems_HTPM}
\left\{\begin{array}{rl}\displaystyle
-\Delta_{y'} (u^{i,k})'+\nabla_{y'}\pi^{i,k}-2N^2{\rm rot}_{y'} w^{i,k}_3=e_i\delta_{1k}&\hbox{ in }Y'_f,\\
\noame
-\Delta_{y'} u^{i,k}_3-2N^2{\rm Rot}_{y'}(w^{i,k})'=0&\hbox{ in }Y'_f,
\\
\noame
{\rm div}_{y'} u^{i,k}=0&\hbox{ in }Y_f,\\
\noame
-R_c\Delta_{y'} (w^{i,k})'+4N^2  (w^{i,k})'-2N^2 \,{\rm rot}_{y'} u^{i,k}_3= \,e_i\delta_{2k}&\hbox{ in }Y'_f,\\
\noame
-R_c\Delta_{y'} w^{i,k}_3+4N^2  w^{i,k}_3-2N^2 \,{\rm rot}_{y'}(u^{i,k})'=0&\hbox{ in }Y'_f,
\\
\noame
u^{i,k}=w^{i,k}=0&\hbox{ in }Y'_s,\\
\noame
\displaystyle\int_{Y'_f}u_3^{i,k}dy'=\int_{Y_f}w_3^{i,k}dy'=0,\\
u^{i,k}(y'), w^{i,k}(y'),\pi^{i,k}(y')\quad Y'-\hbox{periodic}.
\end{array}\right.
\end{equation}

We give the main result concerning the homogenized flow.  

\begin{theorem}\label{them_main_HTPM}
Let $(\hat U, \hat W, \tilde P)\in L^2(\omega;H^1_{\#}(Y')^3)\times  L^2(\omega;H^1_{\#}(Y')^3)\times (L_0^2(\omega)\cap H^1(\omega))$ be the unique weak solution of problem (\ref{hom_system_HTPM}). Then, the extensions $(a_\ep^{-2}\tilde u_\ep,a_\ep^{-1}\tilde w_\ep)$ and  $\tilde P_\ep$ of the solution of problem (\ref{system_2})-(\ref{bc_system_2}) converge weakly to $(\tilde u,\tilde w)$   in  $L^2(\Omega)^3\times L^2(\Omega)^3$ and strongly to $\tilde P$ in $L^2(\omega)$ respectively, with $\tilde u_3=\tilde w_3=0$. Moreover, defining $\widetilde U(x')=\int_0^{1}\tilde u(x',y_3)\,dy_3$ and $\widetilde W(x')=\int_0^{1}\tilde w(x',y_3)\,dy_3$, it holds 
\begin{equation}\label{Darcy_law_u_w_HTPM}
\begin{array}{ll}
\widetilde U'(x')=K^{(1)}_0\left(f'(x')-\nabla_{x'}\tilde P(x')\right)+K^{(2)}_0 g(x'),&\quad \widetilde U_3(x')=0\quad \hbox{ in }\omega,\\
\noame
\widetilde W'(x')=L^{(1)}_0\left(f'(x')-\nabla_{x'}\tilde P(x')\right)+L^{(2)}_0 g(x'),&\quad \widetilde W_3(x')=0\quad \hbox{ in }\omega,
\end{array}
\end{equation}
where $K^{(k)}_0$, $L^{(k)}_0\in \mathbb{R}^{2\times 2}$, $k=1,2$, are matrices with coefficients
$$\left(K^{(k)}_0\right)_{ij}=\int_Y u^{i,k}_j(y)\,dy,\quad \left(L^{(k)}_0\right)_{ij}=\int_Y w^{i,k}_j(y)\,dy,\quad i,j=1,2,$$
where $u^{i,k}$, $w^{i,k}$ are the solutions of the local micropolar problems defined in (\ref{local_problems_HTPM}). 

Here, $\tilde P\in H^1(\omega)\cap L^2_0(\omega)$ is the unique solution of the 2D Darcy equation
\begin{equation}\label{Darcy_law_P_HTPM}
\left\{\begin{array}{l}
{\rm div}_{x'}\left( K^{(1)}_0\left(f'(x')-\nabla_{x'}\tilde P(x')\right)+K^{(2)}_0 g(x')\right)=0\quad \hbox{ in }\omega,\\
\noame
\left( K^{(1)}_0\left(f'(x')-\nabla_{x'}\tilde P(x')\right)+K^{(2)}_0 g(x')\right)\cdot n=0\quad \hbox{ in }\partial\omega.
\end{array}\right.
\end{equation}
\end{theorem}

{\bf Proof.} To eliminate the microscopic variable $y'$ in the effective problem (\ref{hom_system_HTPM}), we proceed as for the critical case by considering the local systems (\ref{local_problems_HTPM}). 

Thanks to the identities for the velocity $\int_{Y'_f} \hat U(x',y')\,dy'=\tilde U(x')$ with $\int_{Y'_f}\hat U_3\,dy'=0$ and the analogous one for the microrotation  given in Lemma \ref{lem_asymp_HTPM}, we deduce that $\widetilde U$ and $\widetilde W$ are given by (\ref{Darcy_law_u_w_HTPM}).

Finally, the divergence condition with respect to the variable $x'$ given in  (\ref{div_x_HTPM}) together with the expression of $\widetilde U'(x')$ gives (\ref{Darcy_law_P_HTPM}),   which has a unique solution since $K^{(1)}_0$ is positive definite and then the whole sequence converges, see Part III - Theorem 2.5.2  in \cite{Luka}.\par\hfill$\square$

\begin{remark} We observe that when $N$ is identically zero, taking into account the linear momentum equations from (\ref{hom_system_HTPM}), we can deduce that the Darcy equation (\ref{Darcy_law_P_HTPM})  agrees with the ones obtained in \cite{AnguianoGrau2, Fabricius} in the case {\bf HTPM}.
\end{remark}

\section{The very thin porous medium (VTPM)}\label{VTPM}
It corresponds to the case when the cylinder height is much smaller than the interspatial distance, i.e. $a_\ep\gg \ep$ which is equivalent to $\lambda=+\infty$.

Next, we give some compactness results about the behavior of the extended sequences $(\tilde u_\ep,\tilde w_\ep,\tilde P_\ep)$ and the unfolding functions $(\hat u_\ep,\hat w_\ep, \hat P_\ep)$ satisfying the {\it a priori} estimates given in Lemmas \ref{lemma_estimates} and \ref{lemma_est_P}, and Lemma \ref{estimates_hat}, respectively.
\begin{lemma}\label{lem_asymp_VTPM}
For a subsequence of $\ep$ still denoted by $\ep$, there exist the following functions:
\begin{itemize}
\item[i)] (Velocity) there exist $\tilde u\in H^1_0(0,1;L^2(\omega)^3)$ with $\tilde u_3=0$ and  $\hat u\in H^1_0(0,1;L^2_{\#}(\omega\times Y')^3)$ with  $\hat u=0$ in $\omega\times Y_s$, such that $\int_{Y}\hat u(x',y)dy=\int_0^{1}\tilde u(x',y_3)\,dy_3$ with $\int_{Y}\hat u_3\,dy=0$, $\hat u_3$ independent of $y_3$ and moreover
\begin{eqnarray}
&\displaystyle \ep^{-2}\tilde u_\ep\rightharpoonup (\tilde u',0)\hbox{  in  }H^1(0,1;L^2(\omega)^3),\quad 
\ep^{-2}\hat u_\ep\rightharpoonup \hat u\hbox{  in  }H^1(0,1;L^2(\omega\times Y')^3),&\label{conv_u_VTPM}\\
\noame
&\displaystyle{\rm div}_{x'}\left(\int_0^{1}\tilde u'(x',y_3)\,dy_3\right)=0\hbox{  in  }\omega,\quad 
\left(\int_0^{1}\tilde u'(x',y_3)\,dy_3\right)\cdot n=0\hbox{  in  }\partial\omega\,,&\label{div_x_VTPM}\\
\noame
&\displaystyle{\rm div}_{y'}\hat u'=0\hbox{  in  }\omega\times Y_f,\quad 
{\rm div}_{x'}\left(\int_{Y_f}\hat u'(x',y)\,dy\right)=0\hbox{  in  }\omega,\quad \left(\int_{Y_f}\hat u'(x',y)\,dy\right)\cdot n=0
\hbox{  in  }\partial\omega\,,&\label{div_y_VTPM}
\end{eqnarray}

\item[ii)] (Microrotation)  there exist $\tilde w\in H^1_0(0,1;L^2(\omega)^3)$ with $\tilde w_3=0$ and  $\hat w\in H^1_0(0,1;L^2_{\#}(\omega\times Y')^3)$ with  $\hat w=0$ in $\omega\times Y_s$, such that $\int_{Y_f}\hat w(x',y)dy=\int_0^{1}\tilde w(x',y_3)\,dy_3$ with $\int_{Y_f}\hat w_3\,dy=0$, $\hat w_3$ independent of $y_3$ and moreover
\begin{eqnarray}
&\displaystyle \ep^{-1}\tilde w_\ep\rightharpoonup (\tilde w',0)\hbox{  in  }H^1(0,1;L^2(\omega)^3),\quad 
 \ep^{-1}\hat w_\ep\rightharpoonup \hat w\hbox{  in  }H^1(0,1;L^2(\omega\times Y')^3),&\label{conv_w_VTPM}
\end{eqnarray}
\item[iii)] (Pressure) there exists $\tilde P\in L^2_0(\Omega)$ independent of $y_3$, such that 
\begin{eqnarray}
&\displaystyle\tilde P_\ep\to \tilde P\hbox{  in  }L^2(\omega),\quad \hat P_\ep\to \tilde P\hbox{  in  }L^2(\omega).& \label{conv_P_VTPM}
\end{eqnarray}
\end{itemize}
\end{lemma}
{\bf Proof. } The proof of this result is obtained by arguing similarly to Section 5 in \cite{AnguianoGrau} and Lemma \ref{lem_asymp_crit} of the present paper.\par\hfill$\square$

\begin{theorem}
In the case {\bf VTPM},   the sequence $(\ep^{-2}\hat u_\ep, \ep^{-1}\hat w_\ep)$ converges weakly to $(\hat u,\hat w)$ in $H^1(0,1;L^2(\omega\times Y')^3)\times H^1(0,1;L^2(\omega\times Y')^3)$ and  $\hat P_\ep$ converges strongly to $\tilde P$ in $L^2(\omega)$, where $(\hat u,\hat w, \tilde P)\in H^1_0(\omega;L^2_{\#}(\omega\times Y')^3)\times H^1_0(\omega;L^2_{\#}(\omega\times Y')^3)\times (L^2_0(\omega)\cap H^1(\omega))$ with $\hat u_3=\hat w_3=0$, is the unique solution of the following homogenized system
\begin{equation}\label{hom_system_VTPM}
\left\{\begin{array}{rl}
\displaystyle
-\partial_{y_3} \hat u'+\nabla_{y'} \hat q=2N^2 \,{\rm rot}_{y_3}\hat w'+f'(x')-\nabla_{x'}\tilde P(x')&\hbox{ in }\omega\times Y_f,\\
\noame
-R_c\partial_{y_3}  \hat w'+4N^2 \hat w=2N^2\,{\rm rot}_{y_3}\hat u'+g'(x')&\hbox{ in }\omega\times Y_f,\\
\noame
{\rm div}_{y'}\hat u'=0&\hbox{ in }\omega\times Y_f,\\
\noame
\hat u'=\hat w'=0&\hbox{ in }\omega\times Y_s,\\
\noame
\displaystyle{\rm div}_{x'}\left(\int_{Y}\hat u'(x',y)\,dy\right)=0&\hbox{ in }\omega,\\
\noame
\displaystyle\left(\int_{Y}\hat u'(x',y)\,dy\right)\cdot n=0&\hbox{ on }\partial\omega,\\
\noame
\hat u'(x',y), \hat w'(x',y),  \hat q(x',y')\quad Y'-\hbox{periodic}.&
\end{array}\right.
\end{equation}
\end{theorem}
{\bf Proof. }We choose $\varphi\in \mathcal{D}(\omega;C_{\#}^\infty(Y)^3)$ with ${\rm div}_{y'}\varphi'=0$ in $\omega\times Y$, ${\rm div}_{x'}(\int_Y \varphi'\,dy)=0$ in $\omega$ and  $\varphi_3$ independent of $y_3$ in (\ref{form_var_hat_u}). Taking into account that thanks to ${\rm div}_{y'}\varphi'=0$ in $\omega\times Y$ and $\varphi_3$ independent of $y_3$, we have that 
$${1\over a_\ep}\int_{\omega\times Y} \hat P_\ep {\rm div}_{y'}\varphi'\,dx'dy=0\quad\hbox{and}\quad {1\over \ep}\int_{\omega\times Y} \partial_{y_3}\hat P_\ep \partial_{y_3}\varphi_3\,dx'dy=0\,.$$
Thus,  passing to the limit using the convergences (\ref{conv_u_VTPM}), (\ref{conv_w_VTPM}), (\ref{conv_P_VTPM}),  $\ep/a_\ep\to 0$ and using in the limit that $\tilde P$ does not depend on $y$ and ${\rm div}_{x'}(\int_Y \varphi'\,dy)=0$, we obtain
\begin{equation}\label{form_varl_limit_VTPM_u}
\left\{\begin{array}{l}
\displaystyle\int_{\omega\times Y}D_{y'}\hat u':D_{y'} \varphi'\,dx'dy=2N^2 \int_{\omega\times Y}{\rm rot}_{y_3}\hat w'\cdot \varphi'\,dx'dy+\int_{\omega\times Y}f'\cdot \varphi'\,dx'dy,\\
\noame
\displaystyle\int_{\omega\times Y}\nabla_{y'}\hat u_3:\nabla_{y'} \varphi_3\,dx'dy=0.
\end{array}\right.
\end{equation}

Next, for every $\psi\in \mathcal{D}(\omega;C_\#^\infty(Y)^3)$ with $\psi_3$ independent of $y_3$, we choose $\psi_\ep=a_\ep^{-1}\psi$ in (\ref{form_var_hat_w}) with $g_\ep$ and $R_M$ satisfying (\ref{estim_f_g_cases_sup}) and (\ref{R_M_super}). Then, passing to the limit using convergences (\ref{conv_u_VTPM}) and (\ref{conv_w_VTPM}), we get
\begin{equation}\label{form_var_limit_VTPM_w}
\left\{\begin{array}{l}\displaystyle R_c\int_{\omega\times Y_f}\partial_{y_3}\hat w':\partial_{y_3}\psi'\,dx'dy+4N^2\int_{\omega\times Y_f}\hat w'\cdot \psi'\,dx'dy={2N^2}\int_{\omega\times Y_f}{\rm rot}_{y_3}\hat u'\cdot \psi'\,dx'dy+\int_{\omega\times Y_f}g'(x')\cdot \psi'\,dx'dy\,,\\
\noame
\displaystyle R_c\int_{\omega\times Y_f}\partial_{y_3}\hat w_3:\partial_{y_3}\psi_3\,dx'dy+4N^2\int_{\omega\times Y_f}\hat w_3\cdot \psi_3\,dx'dy=0.
\end{array}\right.
\end{equation}
The second equations of (\ref{form_varl_limit_VTPM_u}) and (\ref{form_var_limit_VTPM_w}) together to the boundary conditions imply  $\hat u_3=\hat w_3=0$. By density, we can deduce that this variational formulation is equivalent to problem (\ref{hom_system_VTPM}).
\par\hfill$\square$

Let us define the local problems which are useful to eliminate the variable $y$ of the previous homogenized problem
and then obtain a Darcy equation for $\tilde P$. We define $\Phi$ and $\Psi$  by 
\begin{equation}\label{Phi}\Phi(N,R_c)={1\over 12}+{R_c\over 4(1-N^2)}-{1\over 4}\sqrt{{N^2 R_c\over 1-N^2}}\coth\left(N\sqrt{{1-N^2\over R_c}}\right)\,,
\end{equation}
\begin{equation}\label{Psi}
\Psi(N,R_c)={\tanh\left(N\sqrt{{1-N^2\over R_c}}\right)\over {1-{N}\sqrt{{1-N^2\over R_c}}\tanh\left(N\sqrt{{1-N^2\over R_c}}\right)}}\,,
\end{equation}
and for every $k=1,2$, we consider the following 2D local micropolar Darcy problems for $\pi^{i,k}$ as follows
\begin{equation}\label{local_problems_VTPM_Darcy}
\left\{\begin{array}{rl}
\displaystyle
-{\rm div}_{y'}\left({1\over 1-N^2}\Phi(N,R_c)\left(\nabla_{y'}\pi^{i,k}(y')+e_i\delta_{1k}\right)\right)=0& \hbox{  in  }Y'_f\,,\\
\noame
\left({1\over 1-N^2}\Phi(N,R_c)\left(\nabla_{y'}\pi^{i,k}(y')+e_i\delta_{1k}\right)\right)\cdot n=0& \hbox{  in  }\partial Y'_s\,.
\end{array}\right.
\end{equation}
It is known that from the positivity of function $\Phi$, problem (\ref{local_problems_VTPM}) has a unique solution for $\pi^{i,k}\in H^1_{\#}(Y')$ (see \cite{BayadaChamGam} for more details).

\begin{theorem}\label{them_main_VTPM}
Let $(\hat u, \hat w, \tilde P)\in  H^1_0(\omega;L^2_{\#}(\omega\times Y')^3)\times H^1_0(\omega;L^2_{\#}(\omega\times Y')^3)\times (L^2_0(\omega)\cap H^1(\omega))$ be the unique weak solution of problem (\ref{hom_system_VTPM}). Then, the extensions $(\ep^{-2}\tilde u_\ep, \ep^{-1}\tilde w_\ep)$ and  $\tilde P_\ep$ of the solution of problem (\ref{system_2})-(\ref{bc_system_2}) converge weakly to $(\tilde u,\tilde w)$   in $H^1(0, 1;L^2(\omega)^3)\times H^1(0,1;L^2(\omega)^3)$ and strongly to $\tilde P$ in $L^2(\omega)$ respectively, with $\tilde u_3=\tilde w_3=0$. Moreover, defining $\widetilde U(x')=\int_0^{1}\tilde u(x',y_3)\,dy_3$ and $\widetilde W(x')=\int_0^{1}\tilde w(x',y_3)\,dy_3$, it holds 
\begin{equation}\label{Darcy_law_u_w_VTPM}
\begin{array}{ll}
\widetilde U'(x')=K^{(1)}_\infty\left(f'(x')-\nabla_{x'}\tilde P(x')\right)+K^{(2)}_\infty g(x'),&\quad \widetilde U_3(x')=0\quad \hbox{ in }\omega,\\
\noame
\widetilde W'(x')=L^{(2)}_\infty g(x'),&\quad \widetilde W_3(x')=0\quad \hbox{ in }\omega,
\end{array}
\end{equation}
where the matrices $K^{(k)}_\infty\in\mathbb{R}^{2\times 2}$, $k=1,2$,  and $L_\infty^{(2)}\in\mathbb{R}^{2\times 2}$ are matrices with coefficients 
\begin{equation}\label{def_K_sub}
\begin{array}{l}
\displaystyle \left(K^{(k)}_\infty\right)_{ij}={1\over 1-N^2}\int_{Y'}\Phi(N,R_c)\left(\partial_{y_i}\pi^{j,k}(y')+\delta_{ij}\delta_{1k}\right)dy',\quad i,j=1,2\,,\\
\noame
\displaystyle \left(L_\infty^{(2)}\right)_{ij}=-{1\over 4N^3}\sqrt{{R_c\over 1-N^2}}\left(\int_{Y'}\Psi(N,R_c)\,dy'\right)\delta_{ij}\,,
\end{array}
\end{equation}
with  $\Phi$ and $\Psi$  given by (\ref{Phi}) and (\ref{Psi}) respectively, and $\pi^{i,k}\in H^1_{\#}(Y')$, $i,k=1,2$, the unique solutions of the local problems (\ref{local_problems_VTPM}).
Here, $\tilde P\in H^1(\omega)\cap L^2_0(\omega)$ is the unique solution of the 2D Darcy problem
\begin{equation}\label{Darcy_law_P_VTPM}
\left\{\begin{array}{l}
{\rm div}_{x'}\left( K^{(1)}_\infty\left(f'(x')-\nabla_{x'}\tilde P(x')\right)+K^{(2)}_\infty g(x')\right)=0\quad \hbox{ in }\omega,\\
\noame
\left( K^{(1)}_\infty\left(f'(x')-\nabla_{x'}\tilde P(x')\right)+K^{(2)}_\infty g(x')\right)\cdot n=0\quad \hbox{ in }\partial\omega.
\end{array}\right.
\end{equation}
\end{theorem}
{\bf Proof. } We proceed as in in the proof of Theorem \ref{them_main_crit} in order to obtain (\ref{Darcy_law_u_w_VTPM}).
Thus, by using (\ref{idendity_for_y}) where $(u^{i,k},w^{i,k},\pi^{i,k})\in H^1_{0,\#}(Y_f)^2\times H^1_{0,\#}(Y_f)^2\times L^2_0(Y_f')$, $i,k=1,2$, is the unique solution of
  \begin{equation}\label{local_problems_VTPM}
\left\{\begin{array}{rl}\displaystyle
-\partial_{y_3} u^{i,k}+\nabla_{y'}\pi^{i,k}-2N^2{\rm rot}_{y_3} w^{i,k}= -e_i\delta_{1k}&\hbox{ in }Y_f,\\
\noame
{\rm div}_{y'} u^{i,k}=0&\hbox{ in }Y_f,\\
\noame
-R_c\partial_{y_3} w^{i,k}+4N^2  w^{i,k}-2N^2 \,{\rm rot}_{y_3} u^{i,k}= -e_i\delta_{2k}&\hbox{ in }Y_f,\\
\noame
u^{i,k}=w^{i,k}=0&\hbox{ in }Y_s,\\
\noame
\displaystyle\int_{Y'}u_3^{i,k}(y')dy'=\int_{Y'}w_3^{i,k}(y')dy'=0,\\
\noame
u^{i,k}(y), w^{i,k}(y),\pi^{i,k}(y')\quad Y'-\hbox{periodic},
\end{array}\right.
\end{equation}
then, thanks to the identities $\int_{Y_f} \hat u(x',y)\,dy=\int_0^{1}\tilde u(x',y_3)\,dy_3$ with $\hat u_3=0$ and $\int_{Y_f}\hat w(x',y)\,dy=\int_0^{1}\tilde w(x',y_3)\,dy_3$ with 
$\hat w_3=0$ given in Lemma \ref{lem_asymp_VTPM}, it holds 
\begin{equation}\label{Darcy_law_u_w_VTPM_proof}
\begin{array}{ll}
\widetilde U'(x')=\displaystyle \int_Y\hat u'(x',y)\,dy=-K^{(1)}_\infty\left(\nabla_{x'}\tilde P(x')-f'(x')\right)+K^{(2)}_\infty g'(x'),&\quad \widetilde U_3(x')=\displaystyle \int_{Y}\hat u_3(x',y')\,dy=0\quad \hbox{ in }\omega,\\
\noame
\displaystyle \widetilde W'(x')=\displaystyle \int_Y\hat w'(x',y)\,dy=-L^{(1)}_\infty\left(\nabla_{x'}\tilde P(x')-f'(x')\right)+L^{(2)}_\infty g'(x'),&\quad \widetilde W_3(x')=\displaystyle \int_Y\hat w_3(x',y)\,dy=0\quad \hbox{ in }\omega,
\end{array}
\end{equation}
where $K^{(k)}_\infty$, $L^{(k)}_\infty$, $k=1,2$, are matrices defined by their coefficients
\begin{equation}\label{K_infty_proof} \left(K^{(k)}_\infty\right)_{ij}=-\int_{Y}u^{i,k}_j(y)\,dy,\quad \left(L^{(k)}_\infty\right)_{ij}=-\int_{Y}w^{i,k}_j(y)\,dy\,,\quad i,j=1,2\,.
\end{equation}
Then, by the divergence condition in the variable $x'$ given in (\ref{hom_system_VTPM}), we get the Darcy equation (\ref{Darcy_law_P_VTPM}).

However, we observe that (\ref{local_problems_VTPM}) can be viewed as a system of ordinary differential equations with constant coe\-fficients, with respect to the variable $y_3$ and unkowns functions $y_3\mapsto u^{i,k}_1 (y',y_3),w^{i,k}_2 (y',y_3),   u^{i,k}_2 (y',y_3)$, $w^{i,k}_1 (y',y_3)$, where $y'$ is a parameter, $y'\in Y'$. Thus, we can give explicit expressions for $u^{i,k}$ and $w^{i,k}$ given in terms of $\pi^{i,k}$ as follows (see \cite{BayadaChamGam} and  \cite{BayadaLuc} for more details):
\begin{equation}\label{expression_cell_u_w}\begin{array}{rl}
u^{i,k}(y)=& {1\over 2(1-N^2)}\left[y_3^2-y_3+{N^2\over k}\left(\sinh({ky_3})-(\cosh(ky_3)-1)\coth{\left({k \over 2}\right)}\right)\right]\left(\nabla_{y'}\pi^{i,k}(y')+e_i\delta_{1k}\right)\\
\noame 
& +{1\over N^2}\left[\left({2N^2\over k}\sinh(ky_3)-2y_3\right)A+{2N^2\over k}(\cosh(ky_3)-1)B-y_3\right]\left(e_i\delta_{2k}\right)^{\perp}\,,\\
\\
w^{i,k}(y)=& {1\over 4(1-N^2)}\left[2y_3+\left(\cosh(ky_3)-1-\sinh{(ky_3)}\coth\left({k \over 2}\right)\right)\right]\left(\nabla_{y'}\pi^{i,k}(y')+e_i\delta_{1k}\right)^\perp\\
\noame
& -{1\over 2N^2}\Big[\cosh(ky_3)A+\sinh(ky_3)B\Big]e_i\delta_{2k}\,,\\
\end{array}
\end{equation}
where $k=\sqrt{{4N^2(1-N^2)\over R_c}}$ and $A$, $B$ are given by

$$
\begin{array}{l}  \displaystyle A(y')={\sinh(k)\over -2\sinh(k)+{4N^2\over k}(\cosh(k)-1)},\quad B(y')={-(\cosh(k)-1)\over -2\sinh(k)+{4N^2\over k}(\cosh(k)-1)}\,.
 \end{array}
 $$
Integrating with respect to $y_3$, we get
\begin{equation}\label{property_int_uik}
\begin{array}{l}\displaystyle
\int_0^{1}u^{i,k}(y',y_3)\,dy_3=-{1\over 1-N^2}\Phi(N,R_c)\left(\nabla_{y'}\pi^{i,k}(y')+e_i\delta_{1k}\right)\,,\\
\noame
\displaystyle\int_0^{1}w^{i,k}(y',y_3)\,dy_3={1\over 4N^3}\sqrt{{R_c\over 1-N^2}}\Psi(N,R_c)e_i\delta_{2k}\,,
\end{array}
\end{equation}
with $\Phi$ and $\Psi$ given by (\ref{Phi}) and (\ref{Psi}), and so that $\pi^{i,k}$ satisfies the  Darcy local problem (\ref{local_problems_VTPM_Darcy}).  Using the expressions of $u^{i,k}$ and $w^{i,k}$ together with (\ref{Darcy_law_u_w_VTPM_proof}), (\ref{K_infty_proof}) and (\ref{property_int_uik}), we easily get (\ref{Darcy_law_u_w_VTPM}),  which has a unique solution since $K^{(1)}_\infty$ is positive definite and then the whole sequence converges. Observe that, from the second equation in (\ref{property_int_uik}) with $k=2$, we have $L^{(1)}_0=0$, which ends the proof.  
\par\hfill$\square$

\begin{remark} We observe that then $R_c$ tends to zero, the function  $\Phi$ given by (\ref{Phi})  becomes identical to $1/12$. In that case, when $N$ is identically zero, taking into account the linear momentum equations from (\ref{hom_system_VTPM}), we can deduce that the Darcy equation (\ref{Darcy_law_P_VTPM})  agrees with the ones obtained in \cite{AnguianoGrau2, Fabricius} in the case {\bf VTPM}. 
\end{remark}

\end{document}